\documentclass[11pt]{article}
\usepackage{amsfonts}
\usepackage{mathrsfs}
\usepackage{subfigure}
\usepackage{graphicx}
\usepackage{epsfig}
\usepackage{float}
\usepackage[usenames]{color}

 \def\p{\partial} \def\nb{\nonumber}
\def\Vh0{\stackrel{\circ}{V}_h} \def\to{\rightarrow}
   
\def\Om{\Omega}  \def\om{\omega} 
\newcommand{\q}{\quad}
 
\def\l{\label}    
\def\D{\end{document}}   
   
\def\m{\mbox}

\newcommand{\lc}
{\mathrel{\raise2pt\hbox{${\mathop<\limits_{\raise1pt\hbox
{\mbox{$\sim$}}}}$}}}

\newcommand{\gc}
{\mathrel{\raise2pt\hbox{${\mathop>\limits_{\raise1pt\hbox{\mbox{$\sim$}}}}$}}}

\newcommand{\ec}
{\mathrel{\raise2pt\hbox{${\mathop=\limits_{\raise1pt\hbox{\mbox{$\sim$}}}}$}}}

\def\bb{\begin{equation}}  \def\ee{\end{equation}}

\def\beqn{\begin{eqnarray}}  \def\eqn{\end{eqnarray}}

\def\beqnx{\begin{eqnarray*}} \def\eqnx{\end{eqnarray*}}

\def\bn{\begin{enumerate}} \def\en{\end{enumerate}}
\def\i{\item}
\def\bd{\begin{description}} \def\ed{\end{description}}
\def\bfg{\begin{figure}} \def\efg{\end{figure}}

\newtheorem{example}{Example}[section]

\newtheorem{remark}{Remark}[section]
\newtheorem{algorithm}{Algorithm}[section]

 \def\x{{\bf x}}

\title{Overlapping Domain Decomposition Methods
\\ for Linear Inverse Problems}
\author{
Daijun Jiang\footnote{School of
Mathematics and Statistics, Central China Normal
University, Wuhan 430079, PR China. The
work of this author was supported by China Postdoctoral Science Foundation
(Grant no. 2012M521444) and National Natural Science
Foundation of China (Nos 91130022, 11161130003 and 11101317).}
\and
Hui Feng\footnote{School of Mathematical Sciences, Wuhan University,
Wuhan 430072, China. The work of this author was supported by
National Natural Science Foundation of China (No.
91130022, No. 10971159 and No. 11161130003) and NCET
of China.({\tt
hfeng.math@whu.edu.cn}). }
\and Jun
Zou\footnote{Department of Mathematics, The Chinese University of
Hong Kong, Shatin, Hong Kong. {\tt (zou@math.cuhk.edu.hk}). The work
of this author was substantially supported by Hong Kong RGC grants
(Projects 405110 and 404611). }}

\begin{document}
\maketitle

\begin{abstract}

We shall derive and propose several efficient 
overlapping domain decomposition methods for solving some typical
linear inverse problems, including the identification of the flux,
the source strength and the initial temperature in
second order elliptic and parabolic systems. The methods are iterative, and
computationally very efficient: only local forward and adjoint problems
need to be solved in each subdomain, and the local minimizations have
explicit solutions.
%The ill-posed inverse problems are formulated into the most frequently used and
%reliable output least-squares formulation with Tikhonov
%regularization.
Numerical experiments are provided to demonstrate the
robustness and efficiency of the methods, in particular, the convergences
seem nearly optimal, i.e., they do not deteriorate or deteriorate
only slightly when the mesh size reduces.
\end{abstract}



%\bigskip
%{\bf Correspondence address}: Department of Mathematics, The Chinese
%University of Hong Kong, Shatin, Hong Kong. {\bf Email}:
%zou@math.cuhk.edu.hk.
%
%
%
%
%\newpage
\medskip
\section{Introduction}\l{sec:intro}
\setcounter{equation}{0}
%It is well-known that linear inverse problems are not necessarily easier to solve than
%nonlinear inverse problems. In fact the inverse problems of reconstructing heat source
%strength, initial temperature and fluxes are highly ill-posed and unstable.
Domain decomposition methods (DDMs) have been developed and proved to be one of
the most successful methodologies in the construction of efficient numerical solvers
for solving many boundary value and initial-boundary value problems,
the so-called direct problems;
see \cite{tos04} \cite{xu92} \cite{xu98} and the references therein.
DDMs usually possess two important features for solving a wide class of large-scale direct problems:
first, they are natural parallel
solvers and can be easily implemented in parallel computers;
second, their convergence may be made nearly optimal in the sense
that the resulting convergence rate is nearly independent of
the mesh size.

%DDMs have been developed for solving a wide class of boundary and initial-boundary
%value problems or the so-called direct problems.
However, no much progress has been made
in the construction of efficient DDMs for solving mathematically ill-posed inverse problems, although the inverse 
problems are usually much more challenging and time consuming than their corresponding 
direct problems.
%There are a few investigations on DDMs for some inverse problems.
In \cite{chan03} \cite{tai98}, DDMs were used indirectly for an elliptic identification
problem, where classical iterative optimization algorithms were first applied
for the stabilized minimization system of the identification problem,
then the existing DDMs were introduced for solving the direct problems and their adjoint systems
involved at each iteration.
As the outer global iterations of these methods are based on the classical nonlinear optimization
algorithms, their convergences deteriorate rapidly as the degrees of freedom
of the entire optimization systems increase.
Newton's method was first used in \cite{cai09} for solving the optimality system of the 
stabilized minimization of an elliptic identification problem, then an additive Schwarz 
type preconditioned algorithm was applied to solve the linear system involved at each Newton's iteration. 
As Newton's method requires the evaluations of the Hessian of the corresponding objective functional, 
the approach of \cite{cai09} is applicable only to a very special formulation of the parameter identification problem.  
In this work we shall develop some DDMs for directly solving the stabilized minimization systems
of some typical linear inverse problems so that their convergences do not deteriorate or deteriorate
only mildly
as the entire degrees of freedom of the optimization system grow.
Next, we shall briefly address some major difficulties in the construction of
DDMs for inverse problems directly, then point out the new contributions of this work.

We shall use $q$ and $u(q)$ to represent respectively the parameter function
to be identified and the solution to the forward model system
associated with the parameter $q$, then one may formulate a general inverse problem
formally as the following forward operator equation
\beqnx
u(q) = z^\delta %~~~~~~~{\rm or}~~~~~~ u(q)(\cdot, T) = z^\delta,
\eqnx
where $z^\delta$ is the measured data of the exact solution $u$ in some subregion inside the
physical domain or on part of the boundary, or at the terminal
time $t = T$ when the problem is time-dependent.
And the parameter $\delta$ is used here to emphasize the existence of the noise 
in the measured data.

Inverse problems are usually ill-posed as at least one of the following
three conditions is violated: the existence, uniqueness and stability of 
solutions\,\cite{aster05}\cite{bank89}\cite{engl2000}. Of the three conditions stability is the most frequently
encountered difficulty in numerical solutions of inverse problems. One of the most stable and effective approaches to
solve general ill-posed inverse problems is to transform them into stabilized output
least-squares minimizations with some appropriately selected Tikhonov regularizations, namely to minimize
the following type of functionals over some constrained set $K$:
\bb
J(q) = \|u(q)-z^\delta\|^2_V+\beta N(q)
%~~{\rm or}~~ J(q) = \|u(q)(\cdot, T)-z^\delta\|^2_V +\beta N(q),
\label{add}
\ee
where $V$ is a Hilbert or Banach space over the measurement subregion and is determined
based on the type of measurement data available, $N(q)$ is
the regularization term and $\beta$ is a regularization
parameter to balance between the data fitting and regularization.

One of the major difficulties in the construction of DDMs for solving a nonlinear minimization problem
associated with $J(q)$ lies in the global dependence of the
forward operator $u(q)$ on the parameter $q$: a change of $q$ in a small subregion of
the global domain $\Om$  causes the change of $u$ in the entire $\Om$.
This is generally true no matter if $u(q)$ is linear or nonlinear.
Due to this global dependence, a direct application of the DDM
principle to solve the nonlinear minimization problem of $J(q)$
may not work. To illustrate this point more clearly, we consider a decomposition of the global minimization
of $J(q)$ %over a constrained set $K$
over the entire domain $\Om$ into a set of subproblems that
involve only all sub-minimizations of functionals $J_i(q_i+\tilde q)$ %over a constrained set $K_i$
on the subdomains $\Om_i$, where $q_i$ has support only in $\Om_i$, and $\tilde q$ is the known contribution from
other subdomains, then $J_i(q_i+\tilde q)$ should be of the form
\bb
J_i(q_i+\tilde q) = \|u(q_i + \tilde q)-z^\delta\|^2_
{V_i} +\beta N(q_i + \tilde q).
\label{add25}
\ee
Clearly the sub-minimization of functional $J_i$ in (\ref{add25}) involves the solution
$u(q_i +\tilde q)$,
which still needs to solve the forward problem in the global domain $\Om$ even when
operator $u(q)$ is linear and only the local quantity $q_i$ needs to update.
Hence the direct application of the DDMs does not really reduce
the global computations to the ones in the local subdomains.

In this study, we will derive and propose several efficient overlapping
DDMs for solving some typical linear inverse problems, including the identification of
the source strength, the initial temperature inside a physical domain,
and the fluxes on (inaccessible) part of the boundary of a physical domain
in second order elliptic and parabolic systems.
These inverse problems are
all ill-posed, especially unstable with respect to the change of the noise in the data \cite{bank89}.
The new algorithms will be constructed in a way that meets
the true spirits of DDMs, namely at each iteration only smaller minimizations are
solved on the subdomains of the original global domain, and
their convergence is nearly optimal in the sense that
the number of the iterations required for a specified accuracy grows nearly independent
of (or very slowly on) the refinement of finite element meshes.

The rest of the paper is arranged as follows.
%In Section~\ref{sec:lin}, we demonstrate the linearity of the solutions to system (\ref{qi1})
%and give its adjoint operator.
In Section\,\ref{sec:tikf}, we propose the Tikhonov regularization
for identifying the source strength.
%and derive the usual convergence rates for the regularized solutions.
In Section\,\ref{sec:localf}, the overlapping domain
decomposition methods are first introduced and local minimizations are studied, then
the algorithms are further improved. In Sections \ref{sec:tik} and \ref{sec:initial},
we derive DDMs for the reconstruction of the fluxes on part of the boundary and the initial temperature
inside a physical domain respectively.
 In Section\,\ref{sec:nsource},
numerical experiments are presented for the identification
of   source strength, fluxes and initial temperature  to illustrate the efficiency and
robustness of the proposed algorithms. Some concluding remarks are given in
Section\,\ref{sec:conclusion}.

Throughout the paper, $C$ is often used for a generic constant. We shall use the symbol
$\langle \cdot, \cdot\rangle$ for the general inner product, and write the norms of the spaces
$H^m(\Omega)$, $L^2(\Omega)$,
$H^{\frac{1}{2}}(\Gamma)$ and $L^2(\Gamma)$ (for some $\Gamma\subset \p\Omega$)  respectively
as $\|\cdot\|_{m, \Omega}$, $\|\cdot\|_{\Omega}$, $\|\cdot\|_{1/2, \Gamma}$ and
$\|\cdot\|_{\Gamma}$.

%Now we treat $u(h)$ as an operator of $h$ from $L^2(\Gamma_1)$ to
%$L^2(\Gamma_0)$, then the previous equation becomes:
%\bb \left\{
%\begin{array} {lcl}
%-\nabla\cdot(a(\x)\nabla u(h))+u(h)=f \,\,\m{in} \,\, \Om, \\
%a(\x)\frac{\p u(h)}{\p n}=g(\x) ~~~\m{on}
%~~ \Gamma_0,~~~~~~~~~\\
%a(\x)\frac{\p u(h)}{\p n}=h(\x) ~~\m{on} ~~ \Gamma_1.~~~~~~~~~~
%\end{array}
%\right. \l{qi2} \ee

\section{Domain decomposition algorithms for the reconstruction of source strengths}\l{sec:tikf}
%{Regularization formulation and some auxiliary estimates of the direct problem}
\setcounter{equation}{0}

The major task of this work is to propose some new overlapping DDMs 
for solving three typical linear inverse problems, including
the identification of the source strength, the flux and the initial temperature.
For ease of exposition,
we shall take the inverse problem of identifying the source strength in a diffusion
system as an example to derive and discuss the new DDMs in more detail in this section,
and address the other two inverse problems in sections \ref{sec:tik} and \ref{sec:initial}.
%But , since the reconstruction of the heat source
%strength and initial temperature inside the physical domain is similar.
%We will give the general overlapping DDMs in the case of identifying the fluxes on part of the boundary of
%the physical domain. The numerical results for the identification of the fluxes, source and initial value
%will be addressed in Section~\ref{sec:numerical}, \ref{sec:source} and \ref{sec:initial} respectively.
Let $\Om$ be an open bounded and connected domain in ${\mathbb R}^d~(d\geq1)$,
with a boundary $\p\Om$. Then we consider the following diffusion system
%\footnote{Please change the lower order term to
%$c(\x)u$ to make it look more natural. Otherwise this term with constant
%coefficient $1$ looks a bit too special. I think this change does not
%affect our algorithms and theories. Please make this change for all related equations
%in the remaining of this paper.
%
%In addition, please look at the original content of my LaTex writing for
%this equation, then you should change all the equations in this format.
%Most of your equations look very messy.}
%
\bb \left\{ \begin{array} {rclll}
-\nabla\cdot(a(\x)\nabla u)+c(\x)u&=&f(\x) \,\,&\m{in}& \,\, \Om, \\
u(\x)&=&g(\x) ~~&\m{on}& ~\p\Om
\end{array}
\right. \l{qi1f}
\ee
where $a(\x)$, $c(\x)$ and $g(\x)$ are
all given functions, and $a(\x)\geq a_1>0$, $c(\x)\geq c_1>0$ in $\Omega$.
Suppose that the source strength $f(\x)$ of the model system is unknown in $\Omega$.
Our inverse problem is to recover the source strength distribution $f(\x)$
in $\Omega$  when the measurement data
of $u$, denoted by $z^\delta$, is available in $\Omega$, or in a subregion $\tilde \Omega$ of $\Omega$.
%
%A forward problem of (\ref{qi1}) is to determine $u$ on $\Gamma_0$
%when the fluxes $h(\x)$ on $\Gamma_1$ is given.
For convenience,
we shall write the solution of system (\ref{qi1f}) as $u(f)$ to emphasize
its dependence on the source strength $f(\x)$.
%
%Now we assume that the fluxes $h(\x)$ is
%unknown on $\Gamma_1$, but the data $u$ is
%available on $\Gamma_0$, which leads to the
%{\rm \bf identification\,problem\,for\,the\,fluxes:} Given the measurement data
%$z^\delta$ of $u$ on $\Gamma_0$, we reconstruct the fluxes distribution
%$h(\x)$ on $\Gamma_1$.
This is a well-known mathematically ill-posed problem.
As in (\ref{add}), we formulate it in a mathematically
stabilized minimization system of the form
\beqn \min_{f\in L^2(\Om)}J(f)
&=& \|u(f)-z^\delta\|_{\Om}^2+\beta\|f\|_{\Om}^2\,.
\l{q1f}
\eqn
Indeed we can show that the minimizer of the system
is stable in the sense that it depends continuously on the change of
the noise in the data $z^\delta$ \cite{ito01} \cite{xie05}.

{\bf Linearity of the forward solutions}. %\l{sec:lin}
The forward solution $u(f)$ of the system (\ref{qi1f}) is basically
linear in terms of $f$. It is easy to check directly that
$$
u(\lambda_1f_1+\lambda_2f_2)=\lambda_1u(f_1) + \lambda_2u(f_2)
\q \forall\, f_1, f_2 \in L^2(\Om) ~~\mbox{and} ~~\lambda_1, \lambda_2 \in R
$$
if and only if $g(\x) = 0$.
This leads us to consider the solution $U$ to the following system:
\bb \left\{ \begin{array} {rclll}
-\nabla\cdot(a(\x)\nabla U)+c(\x)U&=&f(\x) \,\,&\m{in}& \,\, \Om, \\
U&=&0 ~~&\m{on}& ~\p\Om\,.
\end{array}
\right. \l{qi8f}
\ee
%\bb \left\{
%\begin{array} {rclll}
%-\nabla\cdot(a(\x)\nabla  U(h))+c(\x)U(h)&=&0 &\m{in} & \Om, \\
%a(\x)\frac{\p U(h)}{\p n}&=&0 &\m{on}
%& \Gamma_0,\\
%a(\x)\frac{\p U(h)}{\p n}&=&h &\m{on} & \Gamma_1.
%\end{array}
%\right. \l{qi8} \ee
We can verify that $u(f_1)-u(f_2) = U(f_1-f_2)$ for any $f_1, f_2 \in L^2(\Om)$,
or equivalently we have
\bb u(f)= U(f)+u(0). \l{qi9f} \ee
From now on we shall view the solution $U(f)$ to (\ref{qi8f}) as a mapping from $L^2(\Om)$ to
$L^2(\Om)$.

{\bf Adjoint operator}. It is easy to verify that operator $U(f)$ is self-adjoint.
In fact, we have by integration by parts for any $\om\in {L^2(\Om)}$ that
\beqn
\langle f,\, U(\om)\rangle_{L^2(\Om)}
&=&\langle -\nabla\cdot(a(\x)\nabla U( f))+
c(\x) U( f),\, U(\om)\rangle_{L^2(\Om)}\nb\\
%&=&\int_\Om a(\x)\nabla U(f)\cdot\nabla U(\om) dx+\int_\Om  c(\x)U(f)U(\om) dx\nb\\
&=&\langle  U( f),\,-\nabla\cdot(a(\x)\nabla  U(\om))
+c(\x) U(\om)\rangle_{L^2(\Om)}
=\langle  U(f),\,\om\rangle_{L^2(\Om)}.
\l{y1f}
\eqn
%where we have used the Green formula and $U(\om)=0$ on $\p\Om$ in the second
%equation and the Green formula and $U(f)=0$ on $\p\Om$ in the third
%equation.
%\bb \langle
%U(h),\,\om\rangle_{\Gamma_0}=\langle
%h,\,U^*(\om)\rangle_{\Gamma_1},\,\,\forall\,\om\in L^2(\Gamma_0).
%\l{y1} \ee
%From (\ref{qi8}) and by making use of Green
%formula, we have
%\beqnx &&\langle
%h,\,U^*(\om)\rangle_{\Gamma_1}=\langle a(\x)\frac{\p U(h)}{\p
%n},\,U^*(\om)\rangle_{\Gamma_1}
%=\langle a(\x)\frac{\p U(h)}{\p n},\,U^*(\om)\rangle_{\p\Om}\nb\\
%&=&\int_\Om a(\x)\nabla U(h)\cdot\nabla U^*(\om) dx+\int_\Om \nabla\cdot(a(\x)\nabla U(h))U^*(\om) dx\nb\\
%&=&\int_\Om a(\x)\nabla U(h)\cdot\nabla U^*(\om) dx+\int_\Om  c(\x)U(h)U^*(\om) dx\nb\\
%&=&-\int_\Om U(h)\nabla\cdot(a(\x)\nabla  U^*(\om)) dx+\int_{\p
%\Om}a(\x)\frac{\p U^*(\om)}{\p n}U(h)ds+
%   \int_\Om  c(\x)U(h)U^*(\om) dx\nb\\
%&=&\int_\Om U(h)(-\nabla\cdot(a(\x)\nabla U^*(\om))+c(\x)U^*(\om)) dx+
%   \langle a(\x)\frac{\p U^*(\om)}{\p n},U(h)\rangle_{\Gamma_0\cup\Gamma_1}.
%\eqnx
%Hence (\ref{y1}) holds if $U^*(\om)$ satisfies
%\bb
%\left\{ \begin{array} {rclll}
%-\nabla\cdot(a(\x)\nabla U^*(\om))+c(\x)U^*(\om)&=&0 &\m{in} & \Om, \\
%a(\x)\frac{\p U^*(\om)}{\p n}&=&\omega &\m{on}
%& \Gamma_0,\\
%a(\x)\frac{\p U^*(\om)}{\p n}&=&0 &\m{on} & \Gamma_1.
%\end{array}
%\right. \l{y3} \ee

\subsection{Overlapping DDMs with explicit local solvers}\l{sec:localf}
\setcounter{figure}{0}

%\footnote{Dear Daijun, may I have some frank suggestions and advices for you.
%I asked you to follow my writing on the source strength case of this section to rewrite
%the flux and initial temperature sections and greatly simplify your previous writings
%which were like a repetition of the section for source strength.
%Indeed you have made simplifications,  but it seems you did not closely follow
%my suggestions and the formats I used to write, and still write in your own style and in your mind.
%Many necessary details that should be added but you did not, such as
%notation changes, inappropriate notations,
%equations and section number citations, and still many typos, algorithms formulations
%are not quite logic, .and so on. Please compare my writings with your original ones carefully
%to find out the differences so that you may know how to improve your writings in the future.
%It is very important that when you are writing, you should be very patient, very careful and concentrate,
%and particularly thinking carefully, no matter how busy you are.
%All such details have taken me too much time and also delayed the completion of
%our paper greatly. I hope you will soon become much more improved in writing and organization,
%as I am quite annoyed by spending too much time on so many details. {\bf Thanks very much for your
%sincere advice for improving my writings. It will be helpful in all my life.
%}}
Using the relation (\ref{qi9f}) we can rewrite the minimization (\ref{q1f}) as
\bb \min_{f\in L^2(\Om)}J(f)
%&=&\min_{h\in L^2(\Gamma_1)}\|u(h)-z^\delta\|_{\Gamma_0}^2+\beta\|h\|_{\Gamma_1}^2\nb\\
= \|U(f)-z_0^\delta\|_{\Om}^2+\beta\|f\|_{\Om}^2\,,
\l{q1lf}
\ee
with $z_0^\delta = z^\delta-u(0)$.
As $U(f)$ is linear, $J(f)$ is convex with respect to $f$. And the minimizers
of (\ref{q1lf}) exist and are unique.

In this section, we shall derive some effective DDMs to solve the optimization system
(\ref{q1lf}). We shall not intend to solve this optimization system on the global domain $\Omega$,
as most existing numerical solvers do.
%which is usually very expensive due to the ill-posedness nature of the original inverse problem.
Instead we plan to construct some DDMs so that the nonlinear system (\ref{q1lf})
can be effectively solved on local subdomains.
To do so, we divide the global domain $\Omega$ into a finite number of
overlapping subdomains $\Omega_1$, $\Omega_2$, ... , $\Omega_l$, where $l$ is a positive
integer.
%may be adjusted by the number of processors in a  parallel computer.
Though our new DDMs  work for a general number of subdomains,
we shall focus all our discussions only on 4 subdomains with a cross-point
for ease of exposition; see Figure \ref{dd1}.
It is well-known that the case of 4 subdomains with a cross-point is a most
representative case of general multiple subdomains \cite{tos04} \cite{xu98}.
\bfg[H]
\centerline{\includegraphics[width=200pt]{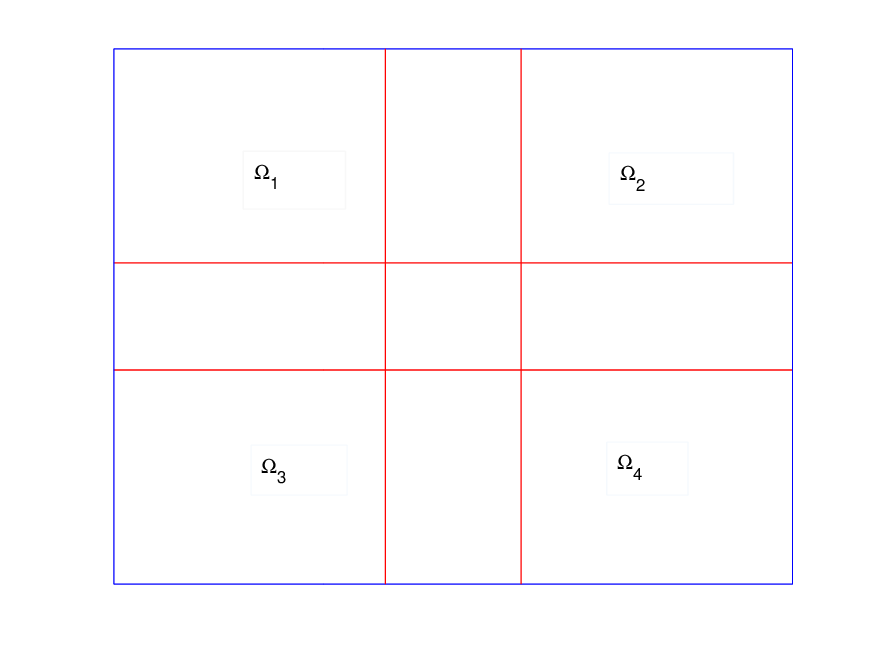}}
\vskip-0.7truecm
\caption{Domain $\Om$ with its $4$ overlapping subdomains $\Om_1, \Om_2, \Om_3, \Om_4$}
\label{dd1}
\efg
Based on the partition of $\Om$ into $4$ overlapping subdomains, we shall often need
a local subspace of $L^2(\Om)$ on each subdomain $\Om_i$:
%\footnote{{\bf we use $V_{f_i}$, $V_{h_i}$ and $V_{\varphi_i}$ to
%distinguish the symbol $V_i$ for three cases. \\
%\mm{This is not a good notation at all: it likes
%the space depending on $f_i$, but actually not. Later we will see big confusion
%when we use $f_i\in V_{f_i}$ or $f_j\in V_{f_j}$. So it is better to change them back
%to $V_i$. Please check the entire paper to change all back to the original notation.}
%} {\bf You are right. I have changed it.}}
\beqnx V_{i}=\Big\{f\in
L^2(\Om); ~{\rm supp}(f)\subset\Om_i\Big\},\, \q
\,i=1,\,2,\,3,\,4. \eqnx
%Noting that when $\p\Om_i\cap\Gamma_1=\emptyset$ for some $i$, then it follows
%that $V_i=\Big\{0\Big\}$.

Next we start to derive some new DD algorithms for solving
the optimization system (\ref{q1lf}). The algorithms are based on
the local optimizations on the subspaces $V_{i}$ associated with
subdomain $\Om_i$.
For some given
%\footnote{the usage of $f_j\in V_{j}$ is extremely confusing: this sounds
%like an element belongs to a space defined by the element, and in your notation one can even
%use $a\in V_a$. Now you may see the notation $V_{f_j}$ is very confusing.
%{\bf You are right. I have changed it.}}
$f_j\in V_{j}$ ($j=1,2,3,4$),
let us consider the following local minimization
over $\Om_i$:
\beqn
\min_{v_i\in V_{i}}J\Big(v_i+ \sum_{j\neq i} f_j\Big)\,.
\l{y4f}
\eqn
Here and in the sequel, we often write $\sum_{j=1,j\neq i}^4$ as $\sum_{j\neq i}$
for simplicity.
By the definition of $J$ in (\ref{q1lf}) we know that each local update $v_i$
in $\Om_i$ still needs to compute the quantity $U(v_i+\sum_{j\neq i} f_j)$,
which involves the solution of the forward system (\ref{qi8f}) in the entire domain $\Om$.
To avoid this, we construct an auxiliary functional
%\footnote{Change all $J_i^s$ to $\tilde J_i^s$
%from here to the first paragraph after equation (\ref{add3f}). {\bf have done.}}
$\tilde J_i^s$ of $J$, called the surrogate functional in \cite{dau04},
by introducing an auxiliary variable $a$.
%\cite{andreas10} \cite{andreas09}.
For a given $a\in V_{i}$ and $f_j\in V_{j}$ ($j=1,2,3,4$),
we define
\bb
\tilde J_i^s(\sum_{j=1}^4f_j,a)=J(\sum_{j=1}^4f_j)+A\,\|f_i-a\|_{\Om}^2-
\|U(f_i-a)\|_{\Om}^2 \l{y5f}
\ee
where $A$ is a positive constant to be selected such that
\bb
A\|f_i-a\|_{\Om}^2-\|U(f_i-a)\|_{\Om}^2\geq
(A-\|U\|^2)\|f_i-a\|_{\Om}^2 \ge 0\,. \label{yan2f}
\ee
This implies for any $f=\sum_{j=1}^4f_j$ that $\tilde J_i^s(f, a)=J(f)$ when $a=f_i$, and
\bb
J(f)=\tilde J_i^s(f,f_i)\leq \tilde J_i^s(f,a)=J(f) + A\|f_i-a\|_{\Om}^2-
\|U(f_i-a)\|_{\Om}^2 \quad \forall\, a\in V_{i}\,.    \label{yan3f}
\ee
%and
%\bb
%%\tilde J_i^s(f,a)-\tilde J_i^s(f,f_i)\geq(A-\|U\|^2)\|f_i-a\|_{\Om}^2.
%\tilde J_i^s(f,a)-J(f)\geq(A-\|U\|^2)\|f_i-a\|_{\Om}^2.
%\label{yan4f}
%\ee
So $\tilde J_i^s(f,a)$ can be viewed as a small perturbation of $J(f)$ when $a$ 
is close to $f_i$. 

Now we shall convert (\ref{y5f}) into a more explicit representation.
Using (\ref{q1lf}), (\ref{y5f}) and the adjoint relation (\ref{y1f})
we can rewrite $\tilde J_i^s$ as follows:
%\footnote{Change all $\sum_{i=1}^4f_i$ into
%$\sum_{j=1}^4f_j$ in the whole paper.}
%\footnote{Please
%write in more details how to get the second and third equalities below. {\bf have done}}
\beqn
\tilde J_i^s(\sum_{j=1}^4f_j,a)
%&=&\| U(f_i)-(z_0^\delta-U(\sum_{j\neq i} f_j))\|^2_{\Om}+\beta\|\sum_{j=1}^4f_j\|_{\Om}^2
%+A\|f_i-a\|_{\Om}^2-\|U(f_i-a)\|_{\Om}^2\nb\\
&=&\| U(f_i)\|_{\Om}^2-2\langle f_i,\,U(z_0^\delta-U(\sum_{j\neq i} f_j)\rangle_{\Om}+
\| z_0^\delta-U(\sum_{j\neq i} f_j)\|^2_{\Om}+\beta\|\sum_{j=1}^4f_j\|_{\Om}^2\nb\\
&&+A\langle f_i,\,f_i-2a\rangle_{\Om}+A\|a\|_{\Om}^2-\| U(f_i)\|_{\Om}^2+
2\langle f_i,\,U(U(a))\rangle_{\Om}-\|U(a)\|_{\Om}^2\nb\\
&=&A\Big\langle f_i,\,f_i-2\Big\{a+\frac{1}{A}U\Big(z_0^\delta-U(\sum_{j\neq i}f_j)-U(a)\Big)\Big\}\Big\rangle_{\Om_i}\nb\\
&&+\beta\|\sum_{j=1}^4f_j\|_{\Om}^2+\| z_0^\delta-U(\sum_{j\neq i} f_j)\|^2_{\Om}
+A\|a\|_{\Om}^2-\|U(a)\|_{\Om}^2\nb\\
&=&A\|f_i-\Big\{a+\frac{1}{A}U\Big(z_0^\delta-U(\sum_{j\neq i}f_j+a)\Big)\Big\}\|_
{\Om_i}^2+\beta\|\sum_{j=1}^4f_j\|_{\Om}^2\nb\\&&+\Big\{\| z_0^\delta-U(\sum_{j\neq i} f_j)\|^2_{\Om}
+A\|a\|_{\Om}^2-\|U(a)\|_{\Om}^2\nb\\&&-A\|a+\frac{1}{A}U\Big(z_0^\delta-U(\sum_{j\neq i}f_j+a)\Big)\|_
{\Om_i}^2\Big\}\,. \l{yyf}
\eqn
%{where we have used the adjoint relation (\ref{y1f}) in the second equation,
%and $f_i\in V_i$ in the third equation.}
%where\footnote{Give explicit formulae for $\varphi$ and $\Phi$ .} $\varphi$ and $\Phi$ are functions of
%$z^\delta,u(0), a$ and $h_{j}$ ($j\neq i$) only.
%\footnote{Give explicit formulae for $\varphi$ and $\Phi$. {\bf have done}}

We can see that
the last two terms in (\ref{yyf}) does not depend on $f_i$, so it will not affect the local minimization
over $\Om_i$
if we drop them in the functional $\tilde J_i^s$. This leads us to consider
the following functional for a given $a\in V_{i}$:
\begin{equation}
\min_{f_i\in V_{i}}\tilde J_i^s(f_i+ \sum_{j\neq i} f_j, a)
=\min_{f_i\in V_{i}}A\|f_i-\tilde z_i \|_{\Om_i}^2+\beta\|\sum_{j=1}^4f_j
\|_{\Om}^2 \l{y6f}
\end{equation}
where $\tilde z_i$ is given by
\bb \l{zif}
\tilde z_i= a+\frac{1}{A}U(z_0^\delta-U(\sum_{j\neq i} f_j+ a)).
\ee
%Additionally in
%\eqnx
Noting that (\ref{y6f}) is a simple quadratic minimization, we can find its exact minimizer $f_i^*$:
%\footnote{{\bf Here, i suggest to write $f_i^*|_{\Omega_i}=\frac{1}{A+\beta}\Big(A\tilde z_i-
%\beta\sum_{j\neq i} f_j\Big)|_{\Omega_i}$, since $f_i^*\in V_{f_i}$, i.e., $f_i^*\in L^2(\Omega)$ and
%$f_i^*=0$ in $\Omega\backslash\Omega_i$.}}
%
\beqn
f_i^*=\frac{1}{A+\beta}\Big(A\tilde z_i-\beta\sum_{j\neq i} f_j\Big)|_{\Omega_i}.
\label{add3f}
\eqn
%The same arguments work for other subproblems $\Om_j$ for $j=2,3,4$.%, i.e.,
%\beqnx
%J_j^s(\sum_{i=1}^4h_i,a)=J(\sum_{i=1}^4h_i)+\|h_j-a\|_{\Gamma_1}^2-
%\|U(h_j-a)\|_{\Gamma_0}^2,
%\eqnx
%with $a,\,h_j\in V_j$ and
%$h_i\in V_i$ for $i\neq j$.
%For the convergence analysis of every subspace minimization, we refer to\footnote{Make it more clear.}
% \cite{tai03} and \cite{tai031}.
%In the next section, we shall study the convergence properties of algorithm~\ref{al:inisequ}.

Clearly, the new functional $\tilde J_i^s$ in (\ref{y6f})
has an obvious advantage over the functional $J$ in (\ref{q1lf})  or (\ref{y4f}):
it is completely local, and the minimization can be solved explicitly
within the subdomain $\Om_i$. However, for the solution of the local minimization
(\ref{y6f})  we need the data $\tilde z_i$ from (\ref{zif}), which involves the evaluations
of $U(\sum_{j\neq i} f_j+ a)$ and $U(z_0^\delta-U(\sum_{j\neq i} f_j+ a))$.
Unfortunately, these two evaluations are both global, and require the solutions of the forward system (\ref{qi8f})
  in the entire domain $\Om$.
This is surely not expected in an efficient DD algorithm.

Next, we shall propose some techniques to get rid of the aforementioned two global evaluations
so that the resulting DD algorithm involves only local minimizations over the local subdomains.
For convenience, we write the boundary of $\Om_i$ inside $\Om$
by $\tilde\Gamma_i$, i.e.,
$\tilde\Gamma_i=\p\Omega_i\cap \Omega$
for $i=1, 2, 3, 4$. Then we introduce a local forward operator
$U_i(f,p)$  associated with the forward problem
(\ref{qi8f}):
\beqn
\left\{ \begin{array} {clclc}
-\nabla\cdot(a(\x)\nabla  U_i( f,\,p))
+c(\x) U_i( f,\,p)
&=& f &\m{in} & \Om_i, \\
 U_i( f,p)&=&0 &\m{on}
& \p\Om\cap\p\Om_i\\
 U_i( f,p)&=&p &\m{on} &\tilde{\Gamma}_i.
\end{array}
\right.
\l{zh1f}
\eqn
%\bb
%\left\{ \begin{array} {clclc}
%-\nabla\cdot(a(\x)\nabla  U_i(h,p))+c(\x)U_i(h,p)&=&0 &\m{in} & \Om_i, \\
%a(\x)\frac{\p U_i(h,p)}{\p n}&=&0 &\m{on}
%& \Gamma_0\cap\p\Om_i,\\
%a(\x)\frac{\p U_i(h,p)}{\p n}&=&h &\m{on}
%& \Gamma_1\cap\p\Om_i,\\
%U_i(h,p)&=&p &\m{on} &\tilde{\Gamma}_i
%%\\
%%U_i(h,p)&=&0 &\m{in} &\Om\setminus\bar{\Om}_i,
%\end{array}
%\right. \l{zh1}
%\ee
%and
%\bb \left\{ \begin{array} {clclc}
%-\nabla\cdot(a(\x)\nabla  U^*_i(\om,q))+c(\x)U^*_i(\om,q)&=&0 &\m{in} & \Om_i, \\
%a(\x)\frac{\p U^*_i(\om,q)}{\p n}&=&0 &\m{on}
%& \Gamma_1\cap\p\Om_i,\\
%a(\x)\frac{\p U^*_i(\om,q)}{\p n}&=&\om &\m{on}
%& \Gamma_0\cap\p\Om_i,\\
%U^*_i(\om,q)&=&q &\m{on} &\tilde{\Gamma}_i.%,\\
%%U^*_i(\om,q)&=&0 &\m{in} &\Om\setminus\bar{\Om}_i.
%\end{array}
%\right. \l{zh2} \ee
Clearly we can split $U_i(f,p)$ as
$U_i(f,p)=U_i(f,0)+U_i(0,p)$, and $U_i(f,0)$ is self-adjoint, i.e.,
%admits
%a very useful relation:
%\footnote{\bf add the proof for $\|U_i(f,0)\|_{\Om_i}
%\leq\|f\|_{\Om_i}$}
%For any $\varphi\in H^1_0(\Om_i)$, the variational system associated
%with the problem (\ref{zh1f}) with $p=0$ is
%\beqnx
%\int_{\Om_i} a(\x)\nabla
%U_i(f,0)\cdot\nabla \varphi dx+\int_{\Om_i} c(\x)U_i(f,0)\varphi dx=
%\int_{\Om_i}f\varphi dx.
%\eqnx
%Then we have the same estimation as in Lemma \ref{lem:condf} that
%$\|U_i(f,0)\|_{\Om_i}\leq\frac{1}{\min\{a_1,c_1\}}\|f\|_{\Om_i}$
\beqn
\langle U_i(f,0),\,\om\rangle_{\Om_i}
&=& \langle f,\,U_i(\om,0)\rangle_{\Om_i} \q \forall\,\om\in
L^2(\Om_i). \l{zh3f}
\eqn
%In fact we can readily deduce by using the integration by parts that\footnote{Write the derivation
%in the same order as the relation (\ref{zh3f}), namely $\langle U_i(f,0),\,\om\rangle_{\Om_i}
%= \cdots = \langle f,\,U_i(\om,0)\rangle_{\Om_i} \q \forall\,\om\in L^2(\Om_i)$. {\bf have done.}}
%\beqnx
%\langle U_i(f,0),\,\om\rangle_{\Om_i}&=&\langle U_i(f,0),\,-\nabla\cdot(a(\x)\nabla  U_i( \om,\,0))
%+c(\x) U_i( \om,\,0)\rangle_{\Om_i}\\
%%&=&\int_{\Om_i} -\nabla\cdot(a(\x)\nabla  U_i( \om,\,0)) U_i(f,0)
%%dx+\int_{\Om_i} c(\x) U_i(f,0))U_i(\om,0) dx\\
%&=&\int_{\Om_i} a(\x)\nabla U_i(f,0)\cdot\nabla U_i(\om,0)
%d\x+\int_{\Om_i} c(\x) U_i(f,0))U_i(\om,0) d\x\\
%&=&\langle -\nabla\cdot(a(\x)\nabla  U_i( f,\,0))
%+c(\x) U_i( f,\,0),\,U_i(\om,0)\rangle_{\Om_i}\nb\\
%&=&\langle f,\,U_i(\om,0)\rangle_{\Om_i}\,.
%\eqnx
%\beqnx
%\langle
%f,\,U_i(\om,0)\rangle_{\Om_i}&=&\langle -\nabla\cdot(a(\x)\nabla  U_i( f,\,0))
%+c(\x) U_i( f,\,0),\,U_i(\om,0)\rangle_{\Om_i}\nb\\
%&=&\int_{\Om_i} a(\x)\nabla U_i(f,0)\cdot\nabla U_i(\om,0)
%dx+\int_{\Om_i} c(\x) U_i(f,0))U_i(\om,0) dx\\
%&=&\int_{\Om_i} -\nabla\cdot(a(\x)\nabla  U_i( \om,\,0)) U_i(f,0)
%dx+\int_{\Om_i} c(\x) U_i(f,0))U_i(\om,0) dx\\
%&=&\langle-\nabla\cdot(a(\x)\nabla  U_i( \om,\,0))
%+c(\x) U_i( \om,\,0),\,U_i(f,0)\rangle_{\Om_i}\\
%&=&\langle
%\om,U_i(f,0)\rangle_{\Om_i}\,.
%\eqnx
%where we have used the Green formula and $U_i(\om,0)=0$ on $\p\Om_i$ in the second
%equation and the Green formula and $U_i(f,0)=0$ on $\p\Om_i$ in the third
%equation.

Using the local operators $U_i(f, p)$ in (\ref{zh1f}),
we introduce the following local functional for $f_j\in V_{j}$, $j=1,2,3,4$:
\beqnx
J_i(\sum_{j=1}^4f_j,p)&=&\|U_i(\sum_{j=1}^4f_j,p)-z_0^\delta\|_{\Om_i}^2
+ \beta\|\sum_{j=1}^4f_j\|_{\Om_i}^2\,,
\eqnx
%\beqnx \left\{ \begin{array} {rcl}
%z_i^\delta=z^\delta \,\,\m{on} \,\, \p\Om_i\cap\Gamma_0, \\
%z_i^\delta=0\,\,\m{on}\,\,\Gamma_0\setminus\p\Om_i,
%\end{array}
%\right.
%\eqnx
and its surrogate functional
%\footnote{Change all $\tilde J_i^s$ to $J_i^s$
%and $\tilde z_i$ to $z_i$
%from here to the end of section 2. {\bf have done.}}
$J_i^s$ for any given $a\in V_{i}$:
$$
J_i^s(\sum_{j=1}^4f_j,p,a)
=J_i(\sum_{j=1}^4f_j,p)+A\|f_i-a\|_{\Om_i}^2-
\|U_i(f_i-a,0)\|_{\Om_i}^2\,.
$$
%\footnote{Write more details to derive
%the second equality. {\bf have done}}
%$i=1,2,3,4$, and take $i=1$ for example:
Using the important fact that $U_i(\sum_{j=1}^4f_j,p)=U_i(\sum_{j\neq i}
f_j,p)+U_i(f_i,0)$ and the adjoint relation (\ref{zh3f}),
we can write
\beqn
J_i^s(\sum_{j=1}^4f_j,p,a)
%&=&J_i(\sum_{j=1}^4f_j,p)+A\|f_i-a\|_{\Om_i}^2-
%\|U_i(f_i-a,0)\|_{\Om_i}^2\nb\\
&=&\| U_i(f_i,0)\|_{\Om_i}^2-2\langle f_i,\,U_i\Big(z_0^\delta-U_i(\sum_{j\neq i}
f_j,p),0\Big)\rangle_{\Om_i}\nb\\
&&+\| z_0^\delta-U_i(\sum_{j\neq i} f_j,p)\|^2_{\Om_i}+
\beta\|\sum_{j=1}^4f_j\|_{\Om_i}^2+A\langle f_i,\,f_i-2a\rangle_{\Om_i}+
A\|a\|_{\Om_i}^2\nb\\&&-\| U_i(f_i,0)\|_{\Om_i}^2+
2\langle f_i,\,U_i(U_i(a,0),0)\rangle_{\Om_i}-\|U_i(a,0)\|_{\Om_i}^2\nb\\
&=&A\langle f_i,\,f_i-2\Big\{a+\frac{1}{A}U_i\Big(z_0^\delta-U_i(\sum_{j\neq i}f_j+a,p),0\Big)\Big\}
\rangle_{\Omega_i}+\beta\|\sum_{j=1}^4f_j\|_{\Omega_i}^2\nb\\
&&+\| z_0^\delta-U_i(\sum_{j\neq i} f_j,p)\|^2_{\Om_i}
+A\|a\|_{\Om_i}^2-\|U_i(a,0)\|_{\Om_i}^2\nb\\
&=&A\|f_i-\Big\{a+\frac{1}{A}U_i\Big(z_0^\delta-U_i(\sum_{j\neq i}f_j+a,p),0\Big)\Big\}\|_{\Om_i}^2
+\beta\|\sum_{j=1}^4f_j\|_{\Om_i}^2\nb\\
&&+\Big\{\| z_0^\delta-U_i(\sum_{j\neq i} f_j,p)\|^2_{\Om_i}
+A\|a\|_{\Om_i}^2-\|U_i(a,0)\|_{\Om_i}^2\nb\\&&-
A\|a+\frac{1}{A}U_i\Big(z_0^\delta-U_i(\sum_{j\neq i}f_j+a,p),0\Big)\|_{\Om_i}^2\Big\},
\label{zzf}
\eqn
%A straightforward computation shows that
%\beqn
%\tilde J_i^s(\sum_{j=1}^4h_j,p,a)
%&=&\|h_i-(a+U^*_i\Big(z^\delta-u(0)-U_i(\sum_{j=1,j\neq i}^4h_j+a,p),0\Big))\|_{\Gamma_1}^2\nb\\
%&&+\beta\|\sum_{j=1}^4h_j\|_{\Gamma_1\cap\p\Om_i}^2+\Psi_i,
%\label{zz}
%\eqn
%\footnote{Give an explicit formula. {\bf have done}}
We can easily see that
the last term above does not depend on $f_i$, so it will not affect the local minimization
over $\Om_i$
if we drop them in the functional $J_i^s$. This leads us to consider
the following functional for a given $a\in V_{i}$:
\beqn
&&\min_{f_i\in V_{i}}
J_i^s(\sum_{j=1}^4f_j,p,a)
= \min_{f_i\in V_{i}}A\|f_i- z_i\|_{\Om_i}^2
+\beta\|\sum_{j=1}^4f_j\|_{\Om_i}^2,
\label{add12f}
\eqn
where $ z_i=a+\frac{1}{A}U_i\Big(z_0^\delta-U_i(\sum_{j\neq i}f_j+a,p),0\Big)$.
(\ref{add12f}) is a simple quadratic minimization, and we can find its exact minimizer
$f_i^*$:
\beqn
f_i^*=\frac{1}{A+\beta}
\Big\{A\,a+U_i\Big(z_0^\delta-U_i(\sum_{j\neq i}f_j+a,p),0\Big)-
\beta\sum_{j\neq i}f_j|_{\Om_i}\Big\}.
\label{jiang1f}
\eqn

We can see from this expression that as long as the inner boundary value $p$ is available,
the minimization (\ref{add12f}) does not involve any global data and is completely local.
Noting that $U(f)|_{\Om_i}=U_i(f,U(f))$ and the definitions of $J_i(f,p)$ and
$J_i^s(f,p,a)$, we can connect
%then for\footnote{$a$ must be in $V_i$, but  $a=\sum_{j=1}^4f_j$
%is not such a function. {\bf I have made a mistake, and here it is for $f=\sum_{j=1}^4f_j$.}}
%$f=\sum_{j=1}^4f_j$ we derive\footnote{What do you want to explain using
%this equivalence in (\ref{zh4f}) ? {\bf To see the relationship between $\tilde J_i(f,p)$,
%$J_i^s(f,p,a)$ and $J(f)$ with substituting $\Om_i$ for $\Om$ in (\ref{q1lf})}.}
$J_i(f,p)$ and $J_i^s(f,p,a)$ with functional $J(f)$ (cf.\,(\ref{q1lf})) restricted in $\Om_i$:
\beqn
&&\|U(\sum_{j=1}^4f_j)-z^\delta_0\|_{\Om_i}^2+\beta\|\sum_{j=1}^4f_j\|_{\Om_i}^2\nb\\
&=&\|U_i(\sum_{j=1}^4f_j,U(\sum_{j=1}^4f_j))-z^\delta_0\|_{\Om_i}^2
+\beta\|\sum_{j=1}^4f_j\|_{\Om_i}^2\nb\\
&=&J_i(\sum_{j=1}^4f_j,U(\sum_{j=1}^4f_j))
=J_i^s(\sum_{j=1}^4f_j,U(\sum_{j=1}^4f_j),f_i). \l{zh4f}
\eqn
%Similarly, we have $J(\sum_{j=1}^4h_j)|_{\bar\Om_i}=\tilde
%J_i^s(\sum_{j=1}^4h_j,U(\sum_{j=1}^4h_j),h_i)$, for $i=2,3,4$.
So using (\ref{add12f}), 
we are now ready to apply the multiplicative or additive Schwarz iteration principle \cite{tos04} \cite{xu98}
to establish two DD algorithms for solving the optimization system
(\ref{q1lf}).
%So in light of
%the idea of Schwarz alternating iteration method,
%we change the sequential algorithm \ref{al:initsequ} and parallel
%algorithm \ref{al:inipara}
% into the following modified algorithms, which are only needed to compute
%$U(h^{(0)})$ once in the whole domain for some initial guess $h^{(0)}$.
%We denote by $p_i^{(0)}$  the restriction of
%$U(h^{(0)})$ on $\tilde \Gamma_i$.
%
For the description of the algorithms, we introduce
%\footnote{\bf introduce the index function $n(x)$.}
an index function for any point $\x\in \Om$:
\bb \l{eq:index}
n(\x)=\Big\{i;\,\x\in \Om_i, ~i\in \{1,2,3,4\}\Big\}\,; \q
|n(\x)|= ~\m{number of elements in} ~n(\x)\,.
\ee

\begin{algorithm}[Multiplicative Schwarz Algorithm (MSA)]\l{al:sequf}
%\footnote{The original arrangements of the algorithm are not so logical.
%So I have made some rearrangements, please read carefully
%to see if it is reasonable or some part is inappropriate. {\bf It is
%quite reasonable and explicit.}}

Choose a tolerance parameter
$\epsilon_1>0$,  an initial value
$f^{(0)}=\sum_{i=1}^4  f_i^{(0)}$ with $ f_i^{(0)}\in V_{i}$ ($i=1,2,3,4$), and solve (\ref{qi8f})  for $U(f^{(0)})$;
set $p_i^{(0)}:= U(f^{(0)})|_{\tilde \Gamma_i}$ and $n:=0$.
\bn
\item Compute $f_i^{(n+1)}\in V_{i}$ sequentially for $i=1$ to $4$ by
\beqn
f_i^{(n+1)}= 
{\rm argmin}_{v_i\in V_{i}}J_i^s(\sum_{j<i}f_j^{(n+1)}+v_i+\sum_{j>i}f_j^{(n)},p_i^{(n)},f_i^{(n)});
\label{add1f}
\eqn
update  $U_i$ in $\Om_i$:
\beqnx
U_i^{(n)}=U_i(\sum_{j\leq i}f_j^{(n+1)}+\sum_{j>i}f_j^{(n)},p_i^{(n)});
\eqnx
update the inner boundary values on $\tilde \Gamma_j$ for $j> i$ if
$\tilde \Gamma_j\in \Om_i$:
\beqnx
p_j^{(n)}=
U_i^{(n)}|_{\tilde \Gamma_j}\,.
\eqnx

\item Compute $f^{(n+1)}=\sum_{i=1}^4 f_i^{(n+1)}$.

\item If $\|f^{(n+1)}-f^{(n)}\|_{\Om}\leq\epsilon_1$, stop the iteration;

otherwise update $U_i$ in subdomain $\Om_i$ ($i=1,2,3,4$):
\beqnx
U_i^{(n+1)}=U_i(f^{(n+1)},p_i^{(n)});
\eqnx
%\footnote{\bf update the inner boundary values by index function.}
{update the inner boundary values on} $\tilde \Gamma_i$ ($i=1,2,3,4$):
\beqnx
p_i^{(n+1)}(\x)=\frac{1}{|n(\x)|}\sum_{j\in n(\x)}U_j^{(n+1)}(\x),\,\
\forall\, \x\in \tilde\Gamma_i.
\eqnx

set $n:=n+1$, go to Step 1.

%and
%$\sum_{i=1}^4\|U_i(h^{(n+1)},p_i^{(n+1)})+u(0)-z^\delta\|_
%{\Gamma_0\cap\p\Om_i}\leq\epsilon_2$, stop.
\en
\end{algorithm}

We can easily see that Algorithm\,\ref{al:sequf} is sequential or multiplicative. 
The next algorithm proposes a parallel version of 
Algorithm\,\ref{al:sequf}.
For this purpose, we introduce a bounded uniform partition of unity
$\{\chi_i\}_{i=1}^4$ such that $\sum_{i=1}^4\chi_i=1$ and
$\|\chi_i\|_{\infty}\leq1$ and ${\rm supp}(\chi_i)\subset \Om_i$.

\begin{algorithm}[Additive Schwarz Algorithm (ASA)]\l{al:paraf}
%\footnote{The original arrangements of the algorithm are not so logical.
%So I have made some rearrangements, please read carefully
%to see if it is reasonable or some part is inappropriate.{\bf It is
%quite reasonable and explicit.}}

Choose a tolerance parameter
$\epsilon_1>0$, a relaxation parameter $\lambda\in (0,1)$,  an initial value
$f^{(0)}=\sum_{i=1}^4  f_i^{(0)}$ with $ f_i^{(0)}\in V_{i}$ ($i=1,2,3,4$), and solve (\ref{qi8f})  for $U(f^{(0)})$;
set $p_i^{(0)}:= U(f^{(0)})|_{\tilde \Gamma_i}$ and $n:=0$.
%
%Choose a tolerance parameter
%$\epsilon_1>0$, an initial value
%$f^{(0)}=\sum_{i=1}^4  \bar f_i^{(0)}$ such that $ \bar f_i^{(0)}\in V_{f_i}$
%for $i=1,2,3,4$ and denote by $p_i^{(0)}$  the restriction of
%$U(f^{(0)})$ on $\tilde \Gamma_i$.  For $n\ge 0$, do the following iteration

\bn
\i Compute $f_i^{(n+1)}\in V_{i}$ in parallel for $i=1,2,3,4$ by
\beqn
f_i^{(n+1)}= {\rm argmin}_{v_i\in V_{i}}
J_i^s(\sum_{j\neq i} f_j^{(n)}+v_i,p_i^{(n)}, f_i^{(n)}).
\label{add10f}
\eqn

\item Compute $f^{(n+1)}=\lambda\sum_{i=1}^4 f_i^{(n+1)}
+(1-\lambda)f^{(n)}$.

\item If $\|f^{(n+1)}-f^{(n)}\|_{\Om}\leq\epsilon_1$, stop the iteration;

otherwise update $U_i$ in subdomains $\Om_i$ ($i=1,2,3,4$):
\beqnx
U_i^{(n+1)}=U_i(f^{(n+1)},p_i^{(n)});
\eqnx
update the inner boundary values on $\tilde \Gamma_i$ ($i=1,2,3,4$):
\beqnx
p_i^{(n+1)}(\x)=\frac{1}{|n(\x)|}\sum_{j\in n(\x)}U_j^{(n+1)}(\x) \q
\forall\, \x\in \tilde\Gamma_i.
\eqnx

set $f_i^{(n+1)}:=\chi_if^{(n+1)}$, and $n:=n+1$, go to Step 1.
%Step 4. If $\|f^{(n+1)}-f^{(n)}\|_
%{\Om}\leq\epsilon_1$, stop;
%%and
%%$\sum_{i=1}^4\|U_i(h^{(n+1)},p_i^{(n+1)})+u(0)-z^\delta\|_
%%{\Gamma_0\cap\p\Om_i}\leq\epsilon_2$, stop;
%otherwise update
%$\bar f_i^{(n+1)}=\chi_if^{(n+1)}$.
%
\en
\end{algorithm}

\begin{remark} 
%formulations:
%\footnote{{\bf whether to write $f_i^{(n+1)}|_{\Om_i}$?}}
%\beqnx
%f_i^{(n+1)}|_{\Om_i}&=&\frac{1}{A+\beta}
%\Big\{Af_i^{(n)}|_{\Om_i}+U_i\Big(z^\delta_0-U_i(\sum_{j<i}f_j^{(n+1)}+\sum_{j\geq i}f_j^{(n)},p_i^{(n)}),0\Big)\\
%&&-
%\beta(\sum_{j<i}f_j^{(n+1)}+\sum_{j>i}f_j^{(n)})|_{\Om_i}\Big\}.
%\eqnx
The same as for (\ref{jiang1f})  we have explicit expressions for the minimizers
$f_i^{(n+1)}$ in (\ref{add1f}) and (\ref{add10f}).
In our numerical implementations, we will simply take
the partition of unity $\{\chi_i\}_{i=1}^4$ used in Algorithm\,\ref{al:paraf} as follows:
$$
\chi_i(\x)={1}/{|n(\x)|} \q \m{for} \q \x\in \Om_i\,; \q \chi_i(\x)=0 \q \m{for}
\q \x\in\Om\backslash\bar\Om_i.
$$
%
%From (\ref{jiang1f}), we can see that the minimizers
%$f_i^{(n+1)}$ in (\ref{add10f}) have the following explicit formulations:
%\footnote{{\bf whether to write $f_i^{(n+1)}|_{\Om_i}$?}}
%\beqnx
%f_i^{(n+1)}|_{\Om_i}&=&\frac{1}{A+\beta}
%\Big\{A\bar f_i^{(n)}|_{\Om_i}+U_i\Big(z_0^\delta-U_i(\sum_{j=1}^4\bar f_j^{(n)},p_i^{(n)}),0\Big)-
%\beta\sum_{j\neq i}\bar f_j^{(n)}|_{\Om_i}\Big\}\\
%&=&\frac{1}{A+\beta}
%\Big\{A\bar f_i^{(n)}|_{\Om_i}+U_i\Big(z_0^\delta-U_i(f^{(n)},p_i^{(n)}),0\Big)-
%\beta\sum_{j\neq i}\bar f_j^{(n)}|_{\Om_i}\Big\},
%\eqnx
%where we have used the relation $\sum_{j=1}^4\bar f_j^{(n)}=\sum_{j=1}^4\chi_jf^{(n)}
%=f^{(n)}$ in the second equality.
\end{remark}

%\footnote{Move this section before the current section 4.}

\section{Domain decomposition algorithms for flux reconstruction}\l{sec:tik}
\setcounter{equation}{0}

%\footnote{This is a simple copy of the section for the source, and should be much simplified
%and better stated.}

In this section,
we propose a DD algorithm to solve the inverse problem of identifying fluxes on part of the boundary.
%But , since the reconstruction of the heat source
%strength and initial temperature inside the physical domain is similar.
%We will give the general overlapping DDMs in the case of identifying the fluxes on part of the boundary of
%the physical domain. The numerical results for the identification of the fluxes, source and initial value
%will be addressed in Section~\ref{sec:numerical}, \ref{sec:source} and \ref{sec:initial} respectively.
Let $\Om\subset R^d~(d\geq1)$ be an open bounded and connected domain,
with a boundary $\p\Om$, which splits into two
parts, i.e., $\p\Om=\Gamma_0\cup\Gamma_1$. Then we consider the
following elliptic system
%\footnote{Please change the lower order term to
%$c(\x)u$ to make it look more natural. Otherwise this term with constant
%coefficient $1$ looks a bit too special. I think this change does not
%affect our algorithms and theories. Please make this change for all related equations
%in the remaining of this paper.
%
%In addition, please look at the original content of my LaTex writing for
%this equation, then you should change all the equations in this format.
%Most of your equations look very messy.}
%
\bb \left\{ \begin{array} {rclll}
-\nabla\cdot(a(\x)\nabla u)+c(\x) u &=&f(\x) & \m{in} &\Om\,, \\
a(\x)\frac{\p u}{\p n}&=&g(\x) & \m{on} &\Gamma_0\,, \\
a(\x)\frac{\p u}{\p n} &=&h(\x) & \m{on} &\Gamma_1\,,
\end{array}
\right. \l{qi1}
\ee
where $a(\x)$, $c(\x)$, $f(\x)$, $g(\x)$ are
all given functions, and $a(\x)\geq a_1>0$, $c(\x)\geq c_1>0$ in $\Omega$.
Suppose that the flux $h(\x)$ of the model system is unknown on the inaccessible part $\Gamma_1$
of $\p\Om$, our inverse problem is to recover the flux distribution on $\Gamma_1$  when some measurement data
$u^\delta$ of $u$ is available on the accessible part $\Gamma_0$ of $\p\Om$.
%
%A forward problem of (\ref{qi1}) is to determine $u$ on $\Gamma_0$
%when the fluxes $h(\x)$ on $\Gamma_1$ is given.
We shall write the solution of system (\ref{qi1}) as $u(h)$ to emphasize
its dependence on the flux $h(\x)$.
%
%Now we assume that the fluxes $h(\x)$ is
%unknown on $\Gamma_1$, but the data $u$ is
%available on $\Gamma_0$, which leads to the
%{\rm \bf identification\,problem\,for\,the\,fluxes:} Given the measurement data
%$z^\delta$ of $u$ on $\Gamma_0$, we reconstruct the fluxes distribution
%$h(\x)$ on $\Gamma_1$.

As discussed in section\,\ref{sec:tikf}, we formulate the  ill-posed inverse problem of recovering the flux
into a mathematically stabilized minimization system of the form
\beqn \min_{h\in L^2(\Gamma_1)}J(h)
&=& \|u(h)-z^\delta\|_{\Gamma_0}^2+\beta\|h\|_{\Gamma_1}^2\,.
\l{q1}
\eqn
This formulation is stable in the sense that the minimizer $h$  to (\ref{q1}) depends
continuously on the change of the noise in the data $u^\delta$ \cite{xie05}.

Similarly to the discussions in Section \ref{sec:tikf},  we can write the solution $u(h)$
to (\ref{qi1}) as
%The forward solution $u(h)$ of the system (\ref{qi1}) is also basically
%linear in terms of $h$. It is easy to check directly
%that $u(h)$ is linear, namely
%$$
%u(\lambda_1h_1+\lambda_2h_2)=\lambda_1u(h_1) + \lambda_2u(h_2)
%\q \forall\, h_1, h_2 \in L^2(\Gamma_1) ~~\mbox{and} ~~\lambda_1, \lambda_2 \in R
%$$
%if and only if
%$f(\x)= 0$ and $g(\x) = 0$.
%This leads
\bb u(h)= U(h)+u(0)\,, \l{qi9} \ee
where $U(h)$ is the solution to the following system:
\bb \left\{
\begin{array} {rclll}
-\nabla\cdot(a(\x)\nabla  U)+c(\x)U&=&0 &\m{in} & \Om, \\
a(\x)\frac{\p U}{\p n}&=&0 &\m{on}
& \Gamma_0,\\
a(\x)\frac{\p U}{\p n}&=&h &\m{on} & \Gamma_1\,.
\end{array}
\right. \l{qi8}
\ee

{\bf Adjoint operator of $U(h)$}. For any $\omega\in L^2(\Gamma_0)$,
consider the solution $U^*(\omega) \in H^1(\Om)$ to
the following system:
\bb
\left\{ \begin{array} {rclll}
-\nabla\cdot(a(\x)\nabla U^*(\om))+c(\x)U^*(\om)&=&0 &\m{in} & \Om, \\
a(\x)\frac{\p U^*(\om)}{\p n}&=&\omega &\m{on}
& \Gamma_0,\\
a(\x)\frac{\p U^*(\om)}{\p n}&=&0 &\m{on} & \Gamma_1.
\end{array}
\right. \l{y3} \ee
This mapping $U^*: L^2(\Gamma_0) \to
L^2(\Gamma_1)$ is the adjoint operator of $U$, namely, it holds that 
\bb \langle
U(h),\,\om\rangle_{\Gamma_0}=\langle
h,\,U^*(\om)\rangle_{\Gamma_1}\,\q\forall\,\om\in L^2(\Gamma_0).
\l{y1} \ee
This relation follows directly from (\ref{qi8}), (\ref{y3}) and an application of
integration by parts:
%\footnote{I have corrected some inappropriate writing
%in your original derivation. {\bf I am sorry, i have found the mistake in my original
%derivation.}}
\beqnx
&&\langle U(h),\,\om\rangle_{\Gamma_0}=\langle U(h), a(\x)\frac{\p U^*(\om)}{\p n}\rangle_{\Gamma_0}\\
&=&\langle U(h),\,a(\x)\frac{\p U^*(\om)}{\p n}\rangle_{\p\Om}+\int_\Om U(h)(-\nabla\cdot(a(\x)
\nabla U^*(\om))+c(\x)U^*(\om)) d\x\\
%&=&\int_\Om a(\x)\nabla U(h)\cdot\nabla U^*(\om) d\x+\int_\Om  c(\x)U(h)U^*(\om) d\x\nb\\
&=&\int_\Om (-\nabla\cdot(a(\x)\nabla U(h)) + c(\x)U(h)) U^*(\om) d\x
+ \langle a(\x)\frac{\p U(h)}{\p n},\,U^*(\om)\rangle_{\p\Om}\nb\\
&=&\langle a(\x)\frac{\p U(h)}{\p n},\,U^*(\om)\rangle_{\Gamma_1}=
\langle h,\,U^*(\om)\rangle_{\Gamma_1}.
\eqnx
%\beqnx &&
%\langle
%h,\,U^*(\om)\rangle_{\Gamma_1}=\langle a(\x)\frac{\p U(h)}{\p
%n},\,U^*(\om)\rangle_{\Gamma_1}
%=\langle a(\x)\frac{\p U(h)}{\p n},\,U^*(\om)\rangle_{\p\Om}\nb\\
%&=&\int_\Om a(\x)\nabla U(h)\cdot\nabla U^*(\om) d\x+
%\int_\Om \nabla\cdot(a(\x)\nabla U(h))U^*(\om) d\x\nb\\
%&=&\int_\Om a(\x)\nabla U(h)\cdot\nabla U^*(\om) d\x+\int_\Om  c(\x)U(h)U^*(\om) d\x\nb\\
%&=&-\int_\Om U(h)\nabla\cdot(a(\x)\nabla  U^*(\om)) d\x+\int_{\p
%\Om}a(\x)\frac{\p U^*(\om)}{\p n}U(h)ds+
%   \int_\Om  c(\x)U(h)U^*(\om) d\x\nb\\
%&=&\int_\Om U(h)(-\nabla\cdot(a(\x)\nabla U^*(\om))+c(\x)U^*(\om)) d\x+
%   \langle a(\x)\frac{\p U^*(\om)}{\p n},U(h)\rangle_{\Gamma_0\cup\Gamma_1}.
%\eqnx

%Furthermore, by the same convenient scaling as in (\ref{yy1f}), we can also have
%%
%%{\bf A convenient scaling and assumption}.
%%For the sake of the subsequent convergence analysis
%%of our domain decomposition methods, from now on we shall assume that
%%the linear operator $U$: $L^2(\Gamma_1)\to
%%L^2(\Gamma_0)$ satisfies the inequality:
%%%\footnote{Do we have to do such scaling
%%%in our numerical simulations ? {\bf In numerical simulations for identifying
%%%the source term, we make a scaling.}}:
%\bb
%\|U\|_{\mathscr{L}(L^2(\Gamma_1),L^2(\Gamma_0))}<1. \label{yy1}
% \ee

\subsection{DD algorithms with explicit local solvers}\l{sec:local1}
In this subsection, we follow section \ref{sec:localf} to derive
some overlapping domain decomposition method for solving the minimization
in (\ref{q1}). As in section \ref{sec:localf}, $\Om$ is divided into the overlapping
subdomains $\Om_i$ ($i=1,2,3,4$), accordingly the feasible constraint space
$L^2(\Gamma_1)$ can be decomposed into the subspaces
\beqnx V_{i}=\Big\{h\in L^2(\Gamma_1); ~{\rm supp}(h)\subset\p\Om_i\cap\Gamma_1\Big\},\, \q
\,i=1,\,2,\,3,\,4. \eqnx
%Clearly, we have $V_{i}=\{0\}$ if $\p\Om_i\cap\Gamma_1=\emptyset$ for some $i$.

Next we introduce an auxiliary surrogate functional
$\tilde J_i^s$ of $J(h)$ in (\ref{q1}) for any given $a\in V_{i}$ and $h_j\in V_{j}$ ($j=1,2,3,4$):
\bb
\tilde J_i^s(\sum_{j=1}^4h_j,a)=J(\sum_{j=1}^4h_j)+A\|h_i-a\|_{\Gamma_1}^2-
\|U(h_i-a)\|_{\Gamma_0}^2. \l{y5}
\ee
By similar derivations to (\ref{yyf}) but using the adjoint relation (\ref{y1}),
we can rewrite $\tilde J_i^s$ as
%\footnote{Please
%write in more details how to get the second and third equalities below. {\bf have done}}
\beqn
\tilde J_i^s(\sum_{j=1}^4h_j,a)
%&=&\| U(h_i)-(z_0^\delta-U(\sum_{j\neq i} h_j))\|^2_{\Gamma_0}+\beta\|\sum_{j=1}^4h_j\|_{\Gamma_1}^2
%+A\|h_i-a\|_{\Gamma_1}^2-\|U(h_i-a)\|_{\Gamma_0}^2\nb\\
&=&A\|h_i-\Big\{a+\frac{1}{A}U^*\Big(z_0^\delta-U(\sum_{j\neq i}h_j+a)\Big)\Big\}\|_
{\p\Omega_i\cap\Gamma_1}^2+\beta\|\sum_{j=1}^4h_j\|_{\Gamma_1}^2\nb\\&&+\Big\{\| z_0^\delta-U(\sum_{j\neq i} h_j)\|^2_{\Gamma_0}
+A\|a\|_{\Gamma_1}^2-\|U(a)\|_{\Gamma_0}^2\nb\\&&-A\|a+\frac{1}{A}U^*\Big(z_0^\delta-U(\sum_{j\neq i}h_j+a)\Big)\|_
{\p\Omega_i\cap\Gamma_1}^2\Big\},
\l{yy}
\eqn
where $z_0^\delta = z^\delta-u(0)$.
%{where we have used the adjoint relation (\ref{y1}) in the second equation,
%and $h_i\in V_i$ in the third equation.}
%where\footnote{Give explicit formulae for $\varphi$ and $\Phi$ .} $\varphi$ and $\Phi$ are functions of
%$z^\delta,u(0), a$ and $h_{j}$ ($j\neq i$) only.
%\footnote{Give explicit formulae for $\varphi$ and $\Phi$. {\bf have done}}
We can see that the last two terms do not depend on $h_i$, so we can drop
them in the minimization of functional $\tilde J_i^s$. Hence it leads us to the following local minimization
for any $ a\in V_{i}$:
\beqn
\min_{h_i\in V_{i}}\tilde J_i^s(h_i+ \sum_{j\neq i} h_j, a)
=\min_{h_i\in V_{i}}A\|h_i-\tilde z_i \|_{\p\Omega_i\cap\Gamma_1}^2+\beta\|\sum_{j=1}^4h_j
\|_{\Gamma_1}^2 \l{y6}
\eqn
where $\tilde z_i$ is given by
\bb \l{zi}
\tilde z_i= a+\frac{1}{A}U^*(z_0^\delta-U(\sum_{j\neq i} h_j+ a)).
\ee
%Additionally in
%\eqnx
This is a quadratic minimization, so we can find its exact minimizer $h_i^*$:
%\footnote{{\bf whether to write $h_i^*|_{\p\Omega_i\cap\Gamma_1}$?}}
\beqn
h_i^*
=\frac{1}{A+\beta}\Big(A\tilde z_i-\beta\sum_{j\neq i} h_j\Big)\Big|_{\p\Omega_i\cap\Gamma_1}.
\label{add3}
\eqn
%The same arguments work for other subproblems $\Om_j$ for $j=2,3,4$.%, i.e.,
%\beqnx
%J_j^s(\sum_{i=1}^4h_i,a)=J(\sum_{i=1}^4h_i)+\|h_j-a\|_{\Gamma_1}^2-
%\|U(h_j-a)\|_{\Gamma_0}^2,
%\eqnx
%with $a,\,h_j\in V_j$ and
%$h_i\in V_i$ for $i\neq j$.
%For the convergence analysis of every subspace minimization, we refer to\footnote{Make it more clear.}
% \cite{tai03} and \cite{tai031}.
%In the next section, we shall study the convergence properties of algorithm~\ref{al:inisequ}.
%\section{Modified sequential and parallel algorithms}\l{sec:modified}
%%\setcounter{equation}{0}
%In comparison with the functional $J$ in (\ref{q1l})  or (\ref{y4}),
%we can now see that the new functional $J_i^s$ in (\ref{y6})
%has an obvious advantage: it is completely local, namely
%any iterative update of $h_i$ can be done completely
%in the subdomain $\Om_i$. However, for the solution of the local minimization
%(\ref{y6})  we need the data $z_i$ from (\ref{zi}), which involves the evaluations
%of $U(\sum_{j\neq i} h_j+ a)$ and $U^*(z_0^\delta-U(\sum_{j\neq i} h_j+ a))$.
%Clearly these two evaluations are global, and require the solutions of the forward system (\ref{qi8})
%and its adjoint system (\ref{y3}) respectively in the entire domain $\Om_i$.
%This is surely not expected in an efficient DD algorithm.

We observe that the minimization (\ref{y6}) is completely local, and its solution
can be achieved explicitly within the subdomain $\Om_i$. However,
its solution $h_i^*$ needs the data $\tilde z_i$ from (\ref{zi}), which involves two global
solutions of the forward and adjoint systems (\ref{qi8}) and (\ref{y3}), and is clearly
not expected in an efficient DD algorithm.
Next, we propose some techniques to avoid these two global evaluations
so that the resulting DD algorithm involves only local minimizations over the local subdomains.
To do so, we  introduce two local forward and adjoint
operators $U_i(h,p)$ and $U^*_i(\omega,q)$ associated with the global forward
and adjoint systems (\ref{qi8}) and (\ref{y3}):
\bb
\left\{ \begin{array} {clclc}
-\nabla\cdot(a(\x)\nabla  U_i(h,p))+c(\x)U_i(h,p)&=&0 &\m{in} & \Om_i, \\
a(\x)\frac{\p U_i(h,p)}{\p n}&=&0 &\m{on}
& \Gamma_0\cap\p\Om_i,\\
a(\x)\frac{\p U_i(h,p)}{\p n}&=&h &\m{on}
& \Gamma_1\cap\p\Om_i,\\
U_i(h,p)&=&p &\m{on} &\tilde{\Gamma}_i
%\\
%U_i(h,p)&=&0 &\m{in} &\Om\setminus\bar{\Om}_i,
\end{array}
\right. \l{zh1}
\ee
and
\bb \left\{ \begin{array} {clclc}
-\nabla\cdot(a(\x)\nabla  U^*_i(\om,q))+c(\x)U^*_i(\om,q)&=&0 &\m{in} & \Om_i, \\
a(\x)\frac{\p U^*_i(\om,q)}{\p n}&=&\omega &\m{on}
& \Gamma_0\cap\p\Om_i,\\
a(\x)\frac{\p U^*_i(\om,q)}{\p n}&=&0 &\m{on}
& \Gamma_1\cap\p\Om_i,\\
U^*_i(\om,q)&=&q &\m{on} &\tilde{\Gamma}_i.%,\\
%U^*_i(\om,q)&=&0 &\m{in} &\Om\setminus\bar{\Om}_i.
\end{array}
\right. \l{zh2} \ee
%We can directly see that
%\beqnx
%U_i(h,p)=U_i(h,0)+U_i(0,p)~~~\mbox{and}
%~~~U_i^*(\om,q)=U_i^*(\om,0)+U_i^*(0,\om), \eqnx
%%\beqnx
%%U_i^*(\om,q)=U_i^*(\om,0)+U_i^*(0,\om)~~~\mbox{and}
%%~~~U_i^*(\om,U^*(\om))=U^*(\om)i,
%%\eqnx
%and $U_i(h,0)$ and
%$U_i^*(\om,0)$ are linear with respect to $h$ and $\om$ respectively.
%\footnote{\bf add the proof for $\|U_i(h,0)\|_{\Gamma_0\cap\p\Om_i}
%\leq\|h\|_{\Gamma_1\cap\p\Om_i}$}
%{For any $\varphi\in H^1(\Om_i)$
%with $\varphi=0$ on $\tilde\Gamma_i$, the variational system associated
%with the system of (\ref{zh1}) with $p=0$ is
%\beqnx
%\int_{\Om_i} a(\x)\nabla
%U_i(h,0)\cdot\nabla \varphi dx+\int_{\Om_i}c(\x) U_i(h,0)\varphi dx=
%\int_{\Gamma_1\cap\p\Om_i}h\varphi ds.
%\eqnx
%Then we have the similar estimation as in Lemma \ref{lem:cond} that
%$\|U_i(h,0)\|_{\Gamma_0\cap\p\Om_i}\leq\frac{C_i^2}{\min\{a_1,c_1\}}\|h\|_{\Gamma_1\cap\p\Om_i}$,
%where $C_i$ are the constants in the upper bound of the trace theorem
%$\|U_i(h,0)\|_{1/2, \p\Om_i}\leq C_i\|U_i(h,0)\|_{1,\Om_i}$.}

Using  the systems (\ref{zh1}), (\ref{zh2}) and the integration by parts formula, 
we derive the following important relation that will be needed later on:
\beqn
\langle U_i(h,0),\,\om\rangle_{\Gamma_0\cap\p\Om_i}
&=& \langle h,\,U^*_i(\om,0)\rangle_{\Gamma_1\cap\p\Om_i} \q \forall\,\om\in
L^2(\Gamma_0\cap\p\Om_i). \l{zh3}
\eqn
%This follows directly from the systems (\ref{zh1}), (\ref{zh2}) and the integration by parts formula:
%%\footnote{There is some mistake in the derivations, please see the derivations below (\ref{y1}) and
%%correct it. {\bf I have corrected it at the third equality. I am so sorry.}}
%\beqnx
%&&\langle
%U_i(h,0),\,\om\rangle_{\Gamma_0\cap\p\Om_i}=\langle U_i(h,0), a(\x)\frac{\p U^*_i(\om,0)}{\p
%n}\rangle_{\Gamma_0\cap\p\Om_i}\\
%&=&\langle U_i(h,0),a(\x)\frac{\p U^*_i(\om,0)}{\p n}\rangle_{\p\Om_i}+
%\int_{\Om_i} U_i(h,0)\{-\nabla\cdot(a(\x)\nabla
%U^*_i(\om,0))\\
%&&+c(\x)U^*_i(\om,0)\} d\x\\
%&=&\int_{\Om_i} a(\x)\nabla U_i(h,0)\cdot\nabla U^*_i(\om,0) d\x+\int_{\Om_i}c(\x)U_i(h,0)U^*_i(\om,0) d\x\\
%&=&\int_{\Om_i} \{-\nabla\cdot(a(\x)\nabla U_i(h,0)) + c(\x)U_i(h,0)\} U_i^*(\om,0) d\x
%+ \langle a(\x)\frac{\p U_i(h,0)}{\p n},\,U^*_i(\om,0)\rangle_{\p\Om_i}\nb\\
%&=&
%\langle a(\x) \frac{\p
%U_i(h,0)}{\p n},\,U^*_i(\om,0)\rangle_{\Gamma_1\cap\p\Om_i}
%=\langle h,\,U^*_i(\om,0)\rangle_{\Gamma_1\cap\p\Om_i}.
%\eqnx

By means of the local operators $U_i(h, p)$ in (\ref{zh1}),
we introduce the local functional $J_i(\sum_{j=1}^4h_j,p)$ for $h_j\in V_{j}$ ($j=1,2,3,4$):
\beqnx
J_i(\sum_{j=1}^4h_j,p)&=&\|U_i(\sum_{j=1}^4h_j,p)-z_0^\delta\|_{\Gamma_0\cap\p\Om_i}^2
+ \beta\|\sum_{j=1}^4h_j\|_{\Gamma_1\cap\p\Om_i}^2\,,
\eqnx
%\beqnx \left\{ \begin{array} {rcl}
%z_i^\delta=z^\delta \,\,\m{on} \,\, \p\Om_i\cap\Gamma_0, \\
%z_i^\delta=0\,\,\m{on}\,\,\Gamma_0\setminus\ap\Om_i,
%\end{array}
%\right.
%\eqnx
and a surrogate functional $J_i^s$ for a given $a\in V_{i}$:
%\footnote{Write more details to derive
%the second equality. {\bf have done}}
%$i=1,2,3,4$, and take $i=1$ for example:
\beqn
J_i^s(\sum_{j=1}^4h_j,p,a)
&=&J_i(\sum_{j=1}^4h_j,p)+A\|h_i-a\|_{\Gamma_1\cap\p\Om_i}^2-
\|U_i(h_i-a,0)\|_{\Gamma_0\cap\p\Om_i}^2\,. \l{eq:sf}
\eqn
%&=&\| U_i(h_i,0)\|_{\Gamma_0\cap\p\Om_i}^2-2\langle h_i,\,U_i^*\Big(z_0^\delta-U_i(\sum_{j\neq i}
%h_j,p),0\Big)\rangle_{\Gamma_1\cap\p\Om_i}\nb\\
%&&+\| z_0^\delta-U_i(\sum_{j\neq i} h_j,p)\|^2_{\Gamma_0\cap\p\Om_i}+
%\beta\|\sum_{j=1}^4h_j\|_{\Gamma_1\cap\p\Om_i}^2+\langle h_i,\,h_i-2a\rangle_{\Gamma_1\cap\p\Om_i}+
%\|a\|_{\Gamma_1\cap\p\Om_i}^2\nb\\&&-\| U_i(h_i,0)\|_{\Gamma_0\cap\p\Om_i}^2+
%2\langle h_i,\,U_i^*(U_i(a),0)\rangle_{\Gamma_1\cap\p\Om_i}-\|U_i(a,0)\|_{\Gamma_0\cap\p\Om_i}^2\nb\\
%&=&\langle h_i,\,h_i-2\Big\{a+U_i^*\Big(z_0^\delta-U_i(\sum_{j\neq i}h_j+a,p),0\Big)\Big\}
%\rangle_{\p\Omega_i\cap\Gamma_1}+\beta\|\sum_{i=1}^4h_i\|_{\p\Omega_i\cap\Gamma_1}^2\nb\\
%&&+\| z_0^\delta-U_i(\sum_{j\neq i} h_j,p)\|^2_{\Gamma_0\cap\p\Om_i}
%+\|a\|_{\Gamma_1\cap\p\Om_i}^2-\|U_i(a,0)\|_{\Gamma_0\cap\p\Om_i}^2\nb\\
Using the important fact that $U_i(\sum_{j=1}^4h_j,p)=U_i(\sum_{j\neq i}
h_j,p)+U_i(h_i,0)$ and the adjoint relation (\ref{zh3}),
%\footnote{You should always mention what relations or arguments we have used to derive the following.
%{\bf have done.}},
we can rewrite $J_i^s(\sum_{j=1}^4h_j,p,a)$ as
\beqn
J_i^s(\sum_{j=1}^4h_j,p,a)&=&A\|h_i-\Big\{a+\frac{1}{A}U^*_i\Big(z_0^\delta-U_i(\sum_{j\neq i}h_j+a,p),0\Big)\Big\}\|_{\Gamma_1\cap\p\Om_i}^2
+\beta\|\sum_{j=1}^4h_j\|_{\Gamma_1\cap\p\Om_i}^2\nb\\
&&+\Big\{\| z_0^\delta-U_i(\sum_{j\neq i} h_j,p)\|^2_{\Gamma_0\cap\p\Om_i}
+A\|a\|_{\Gamma_1\cap\p\Om_i}^2-\|U_i(a,0)\|_{\Gamma_0\cap\p\Om_i}^2\nb\\&&-
A\|a+\frac{1}{A}U^*_i\Big(z_0^\delta-U_i(\sum_{j\neq i}h_j+a,p),0\Big)\|_{\Gamma_1\cap\p\Om_i}^2\Big\}\,.
\label{zz}
\eqn
%{where we have used the relation $U_i(\sum_{j=1}^4h_j,p)=U_i(\sum_{j\neq i}
%h_j,p)+U_i(h_i,0)$ and adjoint relation (\ref{zh3}) in the second equation}.
%A straightforward computation shows that
%\beqn
%\tilde J_i^s(\sum_{j=1}^4h_j,p,a)
%&=&\|h_i-(a+U^*_i\Big(z^\delta-u(0)-U_i(\sum_{j=1,j\neq i}^4h_j+a,p),0\Big))\|_{\Gamma_1}^2\nb\\
%&&+\beta\|\sum_{j=1}^4h_j\|_{\Gamma_1\cap\p\Om_i}^2+\Psi_i,
%\label{zz}
%\eqn
%\footnote{Give an explicit formula. {\bf have done}}
As the last term does not depend on $h_i$, we are led to the following quadratic minimization:
\beqn
&&\min_{h_i\in V_{i}}
J_i^s(\sum_{j=1}^4h_j,p,a)
= \min_{h_i\in V_{i}}A\|h_i-z_i\|_{\Gamma_1\cap\p\Om_i}^2
+\beta\|\sum_{j=1}^4h_j\|_{\Gamma_1\cap\p\Om_i}^2
\label{add12}
\eqn
where $z_i$ is given by
$$
z_i=a+\frac{1}{A}U^*_i\Big(z_0^\delta-U_i(\sum_{j\neq i}h_j+a,p),0\Big).
$$
We can easily find the minimizer to the quadratic optimization (\ref{add12}) in
an explicit form:
\beqn
h_i^*=\frac{1}{A+\beta}
\Big\{Aa+U^*_i\Big(z_0^\delta-U_i(\sum_{j\neq i}h_j+a,p),0\Big)-
\beta\sum_{j\neq i}h_j|_{\Gamma_1\cap\p\Om_i}\Big\}.
\label{jiang1}
\eqn

As in Section\,\ref{sec:localf}, we are now ready 
to formulate two new DD algorithms for the minimization system (\ref{q1})  
for identifying the heat flux.
For the description of the DD algorithms, we introduce
%\footnote{\bf introduce the index function $n(x)$.}
an index function for any point $\x\in \Gamma_1$:
$$
n(\x)=\Big\{i;\,\x\in \p\Om_i\cap\Gamma_1, ~i\in \{1,2,3,4\}\Big\}\,; \q
|n(\x)|= ~\m{number of elements in} ~n(\x)\,.
$$

\begin{algorithm}[Multiplicative Schwarz Algorithm (MSA)]\l{al:sequf1}
%\footnote{Please
%formulate this algorithm as I did for Algorithm 2.1. {\bf have done}}

 Choose a tolerance parameter
$\epsilon_1>0$,  an initial value
$h^{(0)}=\sum_{i=1}^4  h_i^{(0)}$ with $ h_i^{(0)}\in V_{i}$ ($i=1,2,3,4$), and solve (\ref{qi8})  for $U(h^{(0)})$;
set $p_i^{(0)}:= U(h^{(0)})|_{\tilde \Gamma_i}$ and $n:=0$.
\bn
\item Compute $h_i^{(n+1)}\in V_{i}$ sequentially for $i=1$ to $4$ by
\beqn
h_i^{(n+1)}= {\rm argmin}_{v_i\in V_{i}}
J_i^s(\sum_{j<i}h_j^{(n+1)}+v_i+\sum_{j>i}h_j^{(n)},p_i^{(n)},h_i^{(n)});
\label{add1}
\eqn
update  $U_i$ in $\Om_i$:
\beqnx
U_i^{(n)}=U_i(\sum_{j\leq i}h_j^{(n+1)}+\sum_{j>i}h_j^{(n)},p_i^{(n)});
\eqnx
update the inner boundary values on $\tilde \Gamma_j$ for $j> i$ if
$\tilde \Gamma_j\in \Om_i$:
\beqnx
p_j^{(n)}=
U_i^{(n)}|_{\tilde \Gamma_j}\,.
\eqnx

\item Compute $h^{(n+1)}=\sum_{i=1}^4 h_i^{(n+1)}$.

\item If $\|h^{(n+1)}-h^{(n)}\|_{\Gamma_1}\leq\epsilon_1$, stop the iteration;

otherwise update $U_i$ in subdomains $\Om_i$ ($i=1,2,3,4$):
\beqnx
U_i^{(n+1)}=U_i(h^{(n+1)},p_i^{(n)});
\eqnx
%\footnote{\bf update the inner boundary values by index function.}
{update the inner boundary values on} $\tilde \Gamma_i$ ($i=1,2,3,4$):
\beqnx
p_i^{(n+1)}(\x)=\frac{1}{|n(\x)|}\sum_{j\in n(\x)}U_j^{(n+1)}(\x),\,\
\forall\, \x\in \tilde\Gamma_i.
\eqnx

set $n:=n+1$, go to Step 1.

%and
%$\sum_{i=1}^4\|U_i(h^{(n+1)},p_i^{(n+1)})+u(0)-z^\delta\|_
%{\Gamma_0\cap\p\Om_i}\leq\epsilon_2$, stop.
\en
\end{algorithm}

The next algorithm proposes a parallel version of Algorithm\,\ref{al:sequf1}.
For this purpose, we introduce a uniform partition of unity
$\{\chi_i\}_{i=1}^4$ such that $\sum_{i=1}^4\chi_i=1$ and
$\|\chi_i\|_{\infty}\leq1$ and ${\rm supp}(\chi_i)\subset \p\Om_i\cap \Gamma_1$.

\begin{algorithm}[Additive Schwarz Algorithm (ASA)]\l{al:paraf1}
%\footnote{Please
%formulate this algorithm as I did for Algorithm 2.2. {\bf have done}}

Choose a tolerance parameter
$\epsilon_1>0$, a relaxation parameter $\lambda\in (0,1)$,  an initial value
$h^{(0)}=\sum_{i=1}^4  h_i^{(0)}$ with $ h_i^{(0)}\in V_{i}$ ($i=1,2,3,4$),
and solve (\ref{qi8})  for $U(h^{(0)})$;
set $p_i^{(0)}:= U(h^{(0)})|_{\tilde \Gamma_i}$ and $n:=0$.
%
%Choose a tolerance parameter
%$\epsilon_1>0$, an initial value
%$f^{(0)}=\sum_{i=1}^4  \bar f_i^{(0)}$ such that $ \bar f_i^{(0)}\in V_{f_i}$
%for $i=1,2,3,4$ and denote by $p_i^{(0)}$  the restriction of
%$U(f^{(0)})$ on $\tilde \Gamma_i$.  For $n\ge 0$, do the following iteration

\bn
\i Compute $h_i^{(n+1)}\in V_{i}$ in parallel for $i=1,2,3,4$ by
\beqn
h_i^{(n+1)}= {\rm argmin}_{v_i\in V_{i}}
J_i^s(\sum_{j\neq i} h_j^{(n)}+v_i,p_i^{(n)}, h_i^{(n)}).
\label{add10}
\eqn

\item Compute $h^{(n+1)}=\lambda\sum_{i=1}^4 h_i^{(n+1)}
+(1-\lambda)h^{(n)}$.

\item If $\|h^{(n+1)}-h^{(n)}\|_{\Gamma_1}\leq\epsilon_1$, stop the iteration;

otherwise update $U_i$ in subdomains $\Om_i$ ($i=1,2,3,4$):
\beqnx
U_i^{(n+1)}=U_i(h^{(n+1)},p_i^{(n)});
\eqnx
update the inner boundary values on $\tilde \Gamma_i$ ($i=1,2,3,4$):
\beqnx
p_i^{(n+1)}(\x)=\frac{1}{|n(\x)|}\sum_{j\in n(\x)}U_j^{(n+1)}(\x) \q
\forall\, \x\in \tilde\Gamma_i.
\eqnx

set $h_i^{(n+1)}:=\chi_ih^{(n+1)}$, and $n:=n+1$, go to Step 1.
%Step 4. If $\|f^{(n+1)}-f^{(n)}\|_
%{\Om}\leq\epsilon_1$, stop;
%%and
%%$\sum_{i=1}^4\|U_i(h^{(n+1)},p_i^{(n+1)})+u(0)-z^\delta\|_
%%{\Gamma_0\cap\p\Om_i}\leq\epsilon_2$, stop;
%otherwise update
%$\bar f_i^{(n+1)}=\chi_if^{(n+1)}$.
%
\en
\end{algorithm}

\begin{remark} The same as for (\ref{jiang1}),  we have explicit expressions
for the minimizers $h_i^{(n+1)}$ in (\ref{add1}) and (\ref{add10}).
%\footnote{{\bf whether to write $f_i^{(n+1)}|_{\Om_i}$?}}
%\beqnx
%f_i^{(n+1)}|_{\Om_i}&=&\frac{1}{A+\beta}
%\Big\{Af_i^{(n)}|_{\Om_i}+U_i\Big(z^\delta_0-U_i(\sum_{j<i}f_j^{(n+1)}+\sum_{j\geq i}f_j^{(n)},p_i^{(n)}),0\Big)\\
%&&-
%\beta(\sum_{j<i}f_j^{(n+1)}+\sum_{j>i}f_j^{(n)})|_{\Om_i}\Big\}.
%\eqnx
In our numerical implementations, we will simply take
the partition of unity $\{\chi_i\}_{i=1}^4$ used in Algorithm\,\ref{al:paraf1} as follows:
$\chi_i(\x)={1}/{|n(\x)|}$ for $\x\in \p\Om_i\cap\Gamma_1$, and $\chi_i(\x)=0$ for
$\x\in\Gamma_1\backslash\p\Om_i$.
\end{remark}

%\begin{remark}
%Similar to (\ref{jiang1}), we can see that the minimizers
%$h_i^{(n+1)}$ in (\ref{add10}) have the following explicit formulations:
%%\footnote{Write down the explicit formulae for $h_i^{(n+1)}$ in (\ref{add1}) and (\ref{add10}).
%%{\bf have done.}}
%\beqnx
%h_i^{(n+1)}|_{\Gamma_1\cap\p\Om_i}&=&\frac{1}{A+\beta}
%\Big\{A\bar h_i^{(n)}|_{\Gamma_1\cap\p\Om_i}+U^*_i\Big(z_0^\delta-U_i(\sum_{j=1}^4\bar h_j^{(n)},p_i^{(n)}),0\Big)-
%\beta\sum_{j\neq i}\bar h_j|_{\Gamma_1\cap\p\Om_i}\Big\}\\
%&=&\frac{1}{A+\beta}
%\Big\{A\bar h_i^{(n)}|_{\Gamma_1\cap\p\Om_i}+U^*_i\Big(z_0^\delta-U_i(h^{(n)},p_i^{(n)}),0\Big)-
%\beta\sum_{j\neq i}\bar h_j|_{\Gamma_1\cap\p\Om_i}\Big\}.
%\eqnx
%\end{remark}

%\footnote{Move this section right after the section for the fluxes.}
\section{Domain decomposition algorithms for the reconstruction of an initial temperature}\l{sec:initial}
\setcounter{equation}{0}
\setcounter{figure}{0}
\setcounter{table}{0}
In this section, we are interested in extending the DD algorithms
proposed in sections\,\ref{sec:tikf} and \ref{sec:tik} for solving the stationary
inverse source and flux problems to a time-dependent inverse problem, 
the identification of the initial temperature in the following heat conduction system:
%\footnote{The source
%function $f(\x)$ can be also a time-dependent function $f(\x, t)$, yes ? If so, please change it
%to $f(\x, t)$. {\bf Yes}}
\bb \left\{ \begin{array} {rcclc}
u_t-\nabla\cdot(a(\x)\nabla u)&=&f(\x,t)\,\,&\m{in}& \,\, \Om\times (0,\,T), \\
u&=&0 ~~&\m{on}&
~\p\Om\times (0,\,T),\\
u(\x,0)&=&\varphi(\x)\,\,&{\rm in}&\,\,\Om.
\end{array}
\right. \l{ya1}
\ee
We assume that some observation data $z^\delta$ of the temperature $u(\x, t)$ are available in $\Om$
or in some small subregion $\omega\subset \Om$, but
with a time history in the range $[T-\sigma, T]$.
%let $\sigma >0$ be a very small time period, and $Q_\sigma=\Om\times[T-\sigma, T]$,
%in which the terminal status observation data $z^\delta$ is available,
The inverse problem of our interest is to recover the initial temperature distribution
$\varphi(\x)$, using the observation data $z^\delta$.
We shall write the solution of system (\ref{ya1}) as $u(\varphi)$ to emphasize
its dependence on the initial temperature $\varphi(\x)$.

As described in Section\,\ref{sec:tikf}, it is easy to verify that
$u(\varphi)= U(\varphi)+u(0)$, where $U(\varphi)$ is
linear with respect to $\varphi$ and satisfies the following system
\bb \left\{ \begin{array} {rclll}
U_t-\nabla\cdot(a(\x)\nabla U)&=&0\,\,&\m{in}& \,\, \Om\times (0,\,T), \\
U&=&0 ~~&\m{on}&
~\p\Om\times (0,\,T),\\
U(\x,0)&=&\varphi(\x)\,\,&{\rm in}&\,\,\Om,
\end{array}
\right. \l{ya1h}
\ee
whose variational formulation  is given by
\beqn
\int_{0}^{T}\int_{\Om}U_t\psi d\x dt+
\int_{0}^{T}\int_{\Om}a(\x)\nabla U\cdot\nabla\psi d\x dt=0 \q
\forall\,\,\psi\in L^2(0,T;H^1_0(\Om))\,.
\label{add14}
\eqn
Let $z_0^\delta=z^\delta-u(0)$, then we can formulate our inverse problem
as the following regularized output least-squares minimization:
\beqn
\min_{\varphi\in L^2(\Om)}J(\varphi)=\min_{\varphi\in
L^2(\Om)}\int_{T-\sigma}^T\| U(\varphi)-z_0^\delta\|_{L^2(\Om)}^2
dt+\beta\|\varphi\|_{L^2(\Om)}^2. \label{ya2}
\eqn

Now we introduce the {adjoint system} of the forward problem (\ref{ya1h}):
%\footnote{The adjoint
%system is always defined also on the entire time interval $[0, T]$. Please check
%if this change causes any problems in our subsequent derivations. If not, change all $t^*$ to T and
%some related modifications. {\bf There is no problem, i have changed $t^*$ to $T$.}}
\bb \left\{ \begin{array} {rcccc}
U^*_t+\nabla\cdot(a(\x)\nabla U^*)&=&0\,\,&\m{in}& \,\, \Om\times (0,\,T), \\
U^*&=&0 \,\,&\m{on}&
~\p\Om\times (0,\,T),\\
U^*(\x,T)&=&\om\,\,&{\rm in}&\,\,\Om,
\end{array}
\right. \l{ya2h}
\ee
which is linear with respect to $\om$. Next we derive a very useful relation: 
\bb
\langle U(\varphi)(\x,t),\,\om\rangle_{L^2(\Om)} =\langle
\varphi,\,U^*(\om)(\x,T-t)\rangle_{L^2(\Om)}\q \forall\,t\in [0,T]\,.
\label{add16h}
\ee
Clearly, this is true for $t=0$ by the initial and terminal conditions in 
(\ref{ya1h}) and (\ref{ya2h}). To verify it for $t\in (0, T]$, we define $U^{*,s}(\om)$ for $s\in (0,T]$:
\bb \left\{ \begin{array} {rcccc}
U^{*,s}_t+\nabla\cdot(a(\x)\nabla U^{*,s})&=&0\,\,&\m{in}& \,\, \Om\times (0,\,s), \\
U^{*,s}&=&0 \,\,&\m{on}&
~\p\Om\times (0,\,s),\\
U^{*,s}(\x,s)&=&\om\,\,&{\rm in}&\,\,\Om.
\end{array}
\right.
\label{che1}
\ee
It is easy to find the following relation, 
\bb
U^*(\om)(\x,T-s)=U^{*,s}(\om)(\x,0),
\label{che2}
\ee
and the variational formulation of (\ref{che1}), 
\beqn
-\int_{0}^{s}\int_{\Om}U^{*,s}_t\psi d\x dt+
\int_{0}^{s}\int_{\Om}a(\x)\nabla U^{*,s}\cdot\nabla\psi d\x dt=0 \q 
\forall\,\psi\in L^2(0,s;H^1_0(\Om)).
\label{add14h}
\eqn
Using $U^{*,s}$ in (\ref{che1}) and its property (\ref{che2}), we see  (\ref{add16h}) 
immediately from the following relation
\beqn
\langle U(\varphi)(\x,s),\,\om\rangle_{L^2(\Om)} =\langle
\varphi,\,U^{*,s}(\om)(\x,0)\rangle_{L^2(\Om)}\,. 
\l{eq:adjoints}
\eqn
To check this relation, we use (\ref{add14}) with the terminal time $T$ replaced by $s$, 
then take $\psi=U^{*,s}$ and integrate by parts
with respect to $t$ to obtain
\beqn
&&\int_{\Om}U(\x,s)U^{*,s}(\x,s) d\x-\int_{\Om}U(\x,0)U^{*,s}(\x,0) d\x-
\int_{0}^{s}\int_{\Om}UU^{*,s}_t d\x dt\nb\\
&+&\int_{0}^{s}\int_{\Om}a(\x)\nabla U\cdot\nabla U^{*,s} d\x dt=0.
\label{add15h}
\eqn
Now the desired relation (\ref{eq:adjoints}) follows readily from the initial and terminal conditions 
in (\ref{ya1h}) and (\ref{che1}) and equation (\ref{add14h}) with $\psi=U$.

Next we shall follow sections\,\ref{sec:tikf} and \ref{sec:tik} to derive
some overlapping domain decomposition method for solving the time-dependent minimization
(\ref{ya2}). As in section \ref{sec:localf}, $\Om$ is divided into the overlapping
subdomains $\Om_i$ ($i=1,2,3,4$), accordingly the feasible constraint space
$L^2(\Om)$ can be decomposed into the following subspaces:
\beqnx
V_{i}=\Big\{\varphi\in
L^2(\Om); ~~{\rm supp}(\varphi)\subset\Om_i\Big\},\, \q
\,i=1,\,2,\,3,\,4.
\eqnx

In order to avoid any global solution of the forward and adjoint systems
(\ref{ya1h}) and (\ref{ya2h}) in our DD algorithms,  we introduce their local variants, namely,
the solutions $U_i(\varphi,p)$ and $U_i^*(\om,p)$ to the following systems:
\beqnx
\left\{ \begin{array} {clclc}
U_i(\varphi,p)(\x,t)_t-\nabla\cdot(a(\x) \nabla U_i(\varphi,p)(\x,t))&=&
0 &\m{in} & \Om_i\times (0,\,T), \\
U_i(\varphi,p)(\x,t)&=&0 &\m{on}
& (\p\Om\cap\p\Om_i)\times (0,\,T),\\
U_i(\varphi,p)(\x,t)&=&p(\x, t) &\m{on}&\tilde\Gamma_i\times (0,\,T),\\
U_i(\varphi,p)(\x,0)&=&\varphi(\x)&{\rm in}&\Om_i
\end{array}
\right.
\eqnx
and
\beqnx
\left\{ \begin{array} {clclc}
U_i^*(\om,p)(\x,t)_t+\nabla\cdot(a(\x) \nabla U_i^*(\om,p)(\x,t))&=&0 &\m{in} & \Om_i\times (0,\,T), \\
U_i^*(\om,p)(\x,t)&=&0 &\m{on}
& (\p\Om\cap\p\Om_i)\times (0,\,T),\\
U_i^*(\om,p)(\x,t)&=&p(\x, t) &\m{on}&\tilde\Gamma_i\times (0,\,T),\\
U_i^*(\om,p)(\x,T)&=&\om(\x)&{\rm in}&\Om_i.
\end{array}
\right.
\eqnx
Noting that $U_i(\varphi,0)=U_i^*(\om,0)=0$ on $\p\Om_i$, we can derive 
as we did for (\ref{add16h}) that
\bb
\langle U_i(\varphi,0)(\cdot,t),\,\om\rangle_{L^2(\Om_i)} =\langle
\varphi,\,U_i^*(\om,0)(\cdot,T-t)\rangle_{L^2(\Om_i)}\,.
\label{add17h}
\ee
Now we can define a local functional $J_i(\sum_{j=1}^4\varphi_j,p)$ for
$\varphi_j\in V_j$ ($j=1,2,3,4$):
%\footnote{We have not discussed any discretization in space and
%time in this paper, so do you have any special purpose to introduce
%the time discretization by the trapezoidal rule here. If not, please
%rewrite all back to the continuous ones without discretizations like
%what we did in the previous two sections. {\bf now it is the continuous case.}}:
\beqnx
&&J_i(\sum_{j=1}^4\varphi_j,p)
=\int_{T-\sigma}^T\|U_i(\sum_{j=1}^4\varphi_j,p)(\cdot,t)-z_0^\delta\|_{L^2(\Om_i)}^2dt+
\beta\|\sum_{j=1}^4\varphi_j\|_{L^2(\Om_i)}^2%\\
%&=&\frac{\sigma}{2}\|U_i(\sum_{j=1}^4\varphi_j,p)(\x,T-\sigma)-z_0^\delta(\x,T-\sigma)\|_{L^2(\Om_i)}^2\\
%&&+ \frac{\sigma}{2}\|U_i(\sum_{j=1}^4\varphi_j,p)(\x,T)-z_0^\delta(\x,T)\|_{L^2(\Om_i)}^2+
%\beta\|\sum_{j=1}^4\varphi_j\|_{L^2(\Om_i)}^2,
\eqnx
and introduce a surrogate functional $J_i^s$ for any
%\footnote{What's about $p$ ? {\bf now it is
%the continuous one and $p$ does not need to discrete.}}
$a\in V_{i}$:
\beqnx
 J_i^s(\sum_{j=1}^4\varphi_j,p,a)&=&J_i(\sum_{j=1}^4\varphi_j,p)
+A\sigma\|\varphi_i-a\|_{L^2(\Om_i)}^2-\int_{T-\sigma}^T\|U_i(\varphi_i-a,0)(\cdot,t)\|
_{L^2(\Om_i)}^2dt.
\eqnx
%\footnote{Follow the part from equation (\ref{eq:sf}) to the end of the last section to rewrite
%the following content, including the two algorithms here for the initial temperature. {\bf have done}}
Using the fact that $U_i(\sum_{j=1}^4\varphi_j,p)=U_i(\sum_{j\neq i}
\varphi_j,p)+U_i(\varphi_i,0)$ and the adjoint relation (\ref{add17h}), we can rewrite
\beqn
&&J_i^s(\sum_{j=1}^4\varphi_j,p,a)\nb\\
&=&\int_{T-\sigma}^T\{\| U_i(\varphi_i,0)(\cdot,t)\|_{\Om_i}^2-
2\langle U_i(\varphi_i,0)(\cdot,t),\,z_0^\delta-U_i(\sum_{j\neq i}
\varphi_j,p)(\cdot,t)\rangle_{\Om_i}\nb\\
&&+\| z_0^\delta-U_i(\sum_{j\neq i} \varphi_j,p)(\cdot,t)\|^2_{\Om_i}\}dt+
\beta\|\sum_{j=1}^4\varphi_j\|_{\Om_i}^2+A\sigma\langle \varphi_i,\,\varphi_i-2a\rangle_{\Om_i}
+A\sigma\|a\|_{\Om_i}^2\nb\\&&-\int_{T-\sigma}^T\{\| U_i(\varphi_i,0)(\cdot,t)\|_{\Om_i}^2-
2\langle U_i(\varphi_i,0)(\cdot,t),\,U_i(a,0)(\cdot,t)\rangle_{\Om_i}
+\|U_i(a,0)(\cdot,t)\|_{\Om_i}^2\}dt\nb\\
&=&A\sigma\langle \varphi_i,\,\varphi_i-2\{ a+\frac{1}{A\sigma}\int_{T-\sigma}^TU_i^*(z_0^\delta-
U_i(\sum_{j\neq i}\varphi_j+a,p)(\cdot,t),0)(\cdot,T-t)dt\}
\rangle_{\Omega_i}\nb\\
&&+\beta\|\sum_{j=1}^4\varphi_j\|_{\Omega_i}^2+
\int_{T-\sigma}^T\{\| z_0^\delta-U_i(\sum_{j\neq i} \varphi_j,p)(\cdot,t)\|^2_{\Om_i}
-\|U_i(a,0)(\cdot,t)\|_{\Om_i}^2\}dt+A\sigma\|a\|_{\Om_i}^2\nb\\
&=&A\sigma\|\varphi_i-\{a+\frac{1}{A\sigma}\int_{T-\sigma}^TU_i^*(z_0^\delta-
U_i(\sum_{j\neq i}\varphi_j+a,p)(\cdot,t),0)(\cdot,T-t)dt\}\|_{\Om_i}^2
\nb\\
&&+\beta\|\sum_{j=1}^4\varphi_j\|_{\Om_i}^2+
\Big\{\int_{T-\sigma}^T\{\| z_0^\delta-U_i(\sum_{j\neq i} \varphi_j,p)(\cdot,t)\|^2_{\Om_i}
-\|U_i(a,0)(\cdot,t)\|_{\Om_i}^2\}dt+A\sigma\|a\|_{\Om_i}^2\nb\\&&-
A\sigma\|a+\frac{1}{A\sigma}\int_{T-\sigma}^TU_i^*(z_0^\delta-
U_i(\sum_{j\neq i}\varphi_j+a,p)(\cdot,t),0)(\cdot,T-t)dt\}\|_{\Om_i}^2\Big\}.
\label{che3}
\eqn
We can easily see that
the last term above does not depend on $\varphi_i$, so it will not affect the local minimization
over $\Om_i$
if we drop the term in the functional $J_i^s$. This leads us to consider
the following functional for a given $a\in V_{i}$:
\beqn
&&\min_{\varphi_i\in V_{i}}
J_i^s(\sum_{j=1}^4\varphi_j,p,a)
= \min_{\varphi_i\in V_{i}}A\sigma\|\varphi_i- z_i\|_{\Om_i}^2
+\beta\|\sum_{j=1}^4\varphi_j\|_{\Om_i}^2,
\label{che4}
\eqn
where $ z_i=a+\frac{1}{A\sigma}\int_{T-\sigma}^TU_i^*(z_0^\delta-
U_i(\sum_{j\neq i}\varphi_j+a,p)(\cdot,t),0)(\cdot,T-t)dt$.
Clearly the minimization (\ref{che4}) is quadratic, so we can find its exact minimizer
$\varphi_i^*$:
\beqn
\varphi_i^*=\frac{1}{A\sigma+\beta}
\Big\{A\sigma a+\int_{T-\sigma}^TU_i^*(z_0^\delta-
U_i(\sum_{j\neq i}\varphi_j+a,p)(\cdot,t),0)(\cdot,T-t)dt-
\beta\sum_{j\neq i}\varphi_j|_{\Om_i}\Big\}.
\label{che5}
\eqn

By means of the local minimizations (\ref{che4}), 
we are now ready to formulate two new DD algorithms for solving the minimization 
(\ref{ya2}) for the reconstruction of the initial temperature.
The same index function $n(\x)$ as in (\ref{eq:index}) is used below for any $\x\in \Om$.
%$$
%n(\x)=\Big\{i;\,\x\in \Om_i, ~i\in \{1,2,3,4\}\Big\}\,; \q
%|n(\x)|= ~\m{number of elements in} ~n(\x)\,.
%$$

\begin{algorithm}[Multiplicative Schwarz Algorithm (MSA)]\l{al:initial}
%%%%%%%%%%%%%

Choose a tolerance parameter
$\epsilon_1>0$,  an initial value
$\varphi^{(0)}=\sum_{i=1}^4  \varphi_i^{(0)}$ with $ \varphi_i^{(0)}\in V_{i}$ ($i=1,2,3,4$),
and solve (\ref{ya1h})  for $U(\varphi^{(0)})$;
set $p_i^{(0)}:= U(\varphi^{(0)})|_{\tilde \Gamma_i}$ and $n:=0$.
\bn
\item Compute $\varphi_i^{(n+1)}\in V_{i}$ sequentially for $i=1$ to $4$ by
\beqn
\varphi_i^{(n+1)}= {\rm argmin}_{v_i\in V_{i}}
J_i^s(\sum_{j<i}\varphi_j^{(n+1)}+v_i+\sum_{j>i}\varphi_j^{(n)},p_i^{(n)},\varphi_i^{(n)});
\label{che6}
\eqn
update  $U_i$ in $\Om_i$:
\beqnx
U_i^{(n)}=U_i(\sum_{j\leq i}\varphi_j^{(n+1)}+\sum_{j>i}\varphi_j^{(n)},p_i^{(n)});
\eqnx
update the inner boundary values on $\tilde \Gamma_j$ for $j> i$ if
$\tilde \Gamma_j\in \Om_i$:
\beqnx
p_j^{(n)}=
U_i^{(n)}|_{\tilde \Gamma_j}\,.
\eqnx

\item Compute $\varphi^{(n+1)}=\sum_{i=1}^4 \varphi_i^{(n+1)}$.

\item If $\|\varphi^{(n+1)}-\varphi^{(n)}\|_{\Om}\leq\epsilon_1$, stop the iteration;

otherwise update $U_i$ in subdomain $\Om_i$ ($i=1,2,3,4$):
\beqnx
U_i^{(n+1)}=U_i(\varphi^{(n+1)},p_i^{(n)});
\eqnx
%\footnote{\bf update the inner boundary values by index function.}
{update the inner boundary values on} $\tilde \Gamma_i$ ($i=1,2,3,4$):
\beqnx
p_i^{(n+1)}(\x)=\frac{1}{|n(\x)|}\sum_{j\in n(\x)}U_j^{(n+1)}(\x),\,\
\forall\, \x\in \tilde\Gamma_i.
\eqnx

set $n:=n+1$, go to Step 1.

%and
%$\sum_{i=1}^4\|U_i(h^{(n+1)},p_i^{(n+1)})+u(0)-z^\delta\|_
%{\Gamma_0\cap\p\Om_i}\leq\epsilon_2$, stop.
\en
\end{algorithm}

The next algorithm is a parallel version of Algorithm\,\ref{al:initial}.
%For this purpose, we introduce a bounded uniform partition of unity
%$\{\chi_i\}_{i=1}^4$ such that $\sum_{i=1}^4\chi_i=1$ and
%$\|\chi_i\|_{\infty}\leq1$ and ${\rm supp}(\chi_i)\subset \Om_i$.

\begin{algorithm}[Additive Schwarz Algorithm (ASA)]\l{al:inip}

Choose a tolerance parameter
$\epsilon_1>0$, a relaxation parameter $\lambda\in (0,1)$,  an initial value
$\varphi^{(0)}=\sum_{i=1}^4  \varphi_i^{(0)}$ with $ \varphi_i^{(0)}\in V_{i}$
($i=1,2,3,4$), and solve (\ref{ya1h})  for $U(\varphi^{(0)})$;
set $p_i^{(0)}:= U(\varphi^{(0)})|_{\tilde \Gamma_i}$ and $n:=0$.

\bn
\i Compute $\varphi_i^{(n+1)}\in V_{i}$ in parallel for $i=1,2,3,4$ by
\beqn
\varphi_i^{(n+1)}= {\rm argmin}_{v_i\in V_{i}}
J_i^s(\sum_{j\neq i} \varphi_j^{(n)}+v_i,p_i^{(n)}, \varphi_i^{(n)}).
\label{che7}
\eqn

\item  Compute $\varphi^{(n+1)}=\lambda\sum_{i=1}^4 \varphi_i^{(n+1)}
+(1-\lambda)\varphi^{(n)}$.

\item If $\|\varphi^{(n+1)}-\varphi^{(n)}\|_{\Om}\leq\epsilon_1$, stop the iteration;

otherwise update $U_i$ in subdomains $\Om_i$ ($i=1,2,3,4$):
\beqnx
U_i^{(n+1)}=U_i(\varphi^{(n+1)},p_i^{(n)});
\eqnx
update the inner boundary values on $\tilde \Gamma_i$ ($i=1,2,3,4$):
\beqnx
p_i^{(n+1)}(\x)=\frac{1}{|n(\x)|}\sum_{j\in n(\x)}U_j^{(n+1)}(\x) \q
\forall\, \x\in \tilde\Gamma_i.
\eqnx

set $\varphi_i^{(n+1)}:=\chi_i\varphi^{(n+1)}$, and $n:=n+1$, go to Step 1.
\en
\end{algorithm}

\begin{remark} 
The same as for (\ref{che5}), we have explicit expressions for the minimizers
$\varphi_i^{(n+1)}$ in (\ref{che6}) and (\ref{che7}).
\end{remark}

\section{Numerical experiments }\l{sec:nsource}
\setcounter{equation}{0}
\setcounter{figure}{0}
\setcounter{table}{0}
In this section, we shall apply the DD algorithms that were proposed in 
the previous Sections~\ref{sec:tikf}-\ref{sec:initial} to identify the source strength 
in the elliptic system (\ref{qi1f}), the heat flux in the system
(\ref{qi1}) and the initial temperature in the parabolic system (\ref{ya1}) respectively.
%We implement the {\bf MSA} and {\bf ASA}  using Matlab.

We choose the domain $\Omega=(0,1)\times(0,2)$ and
decompose it into four overlapping subdomains:
$\Omega_1=(0,4/7)\times(6/7,2)$,
$\Omega_2=(3/7,1)\times(6/7,2)$, $\Omega_3=(0,4/7)\times(0,8/7)$,
$\Omega_4=(3/7,1)\times(0,8/7)$. Then we triangulate domain $\Omega$ into
$N\times M$ small squares of equal size and further divide each
square through its diagonal into two triangles. This results in a finite element triangulation 
of domain $\Om$, which is done in such a way that it is consistent with the 
subdomain decompositions. 
All the elliptic problems involved in DD algorithms are solved by the continuous linear finite
element method, while all the parabolic problems are solved by the continuous linear finite
element method in space and the Crank-Nicolson scheme in time.
%\footnote{We should also mention
%what numerical methods we used to solve all the forward problems involved in
%our numerical experiments. {\bf have done}}

The parameters involved in the DD algorithms are chosen as follows. 
The initial guesses are set to be identically equal to some constants, 
which as we see are rather poor initial guesses for all the test problems.  
We take the parameter $A=1$ and the relaxation parameter $\lambda={1}/{2}$
in all the numerical experiments. The noisy data $z^\delta$ is obtained by adding 
some uniform random noise to the exact data, i.e., 
$z^\delta=u + \delta R\, u$, where $R$ is a uniform random function varying in the
range [-1,1].
%\footnote{The numerical experiments for $\delta=0.05$ will be shown
%later.}
The errors shown in all the tables are the relative $L^2$-norm errors
$\|q^{(k)}-q\|/\|q\|$, where $q$ and $q^{(k)}$ are the exact parameter and
its numerical reconstruction by the DD algorithms, which are terminated  
when the relative $L^2$-norm errors reach $0.1$. The exact parameters and their numerical 
reconstructed profiles will be also presented. 

\medskip
We start two numerical tests for the flux reconstructions in the partial 
boundary $\Gamma_1=\{(x,y); ~x=1,\,\,0\leq y\leq2\}$ in the system
(\ref{qi1}), where we take $g(\x)=0$ on $\Gamma_0$, $f(\x)=0$ and
$a(\x)=c(\x)=1$ in $\Om$.

\begin{example}
We take the exact flux $h=-(y-1)^2+1$ on $\Gamma_1$, and the noise level $\delta=5\%$, 
with the constant initial guess $h^{(0)}=1$.
\end{example}

\begin{figure}
  \centering
  \begin{tabular}{cc}
    \includegraphics[width=7cm]{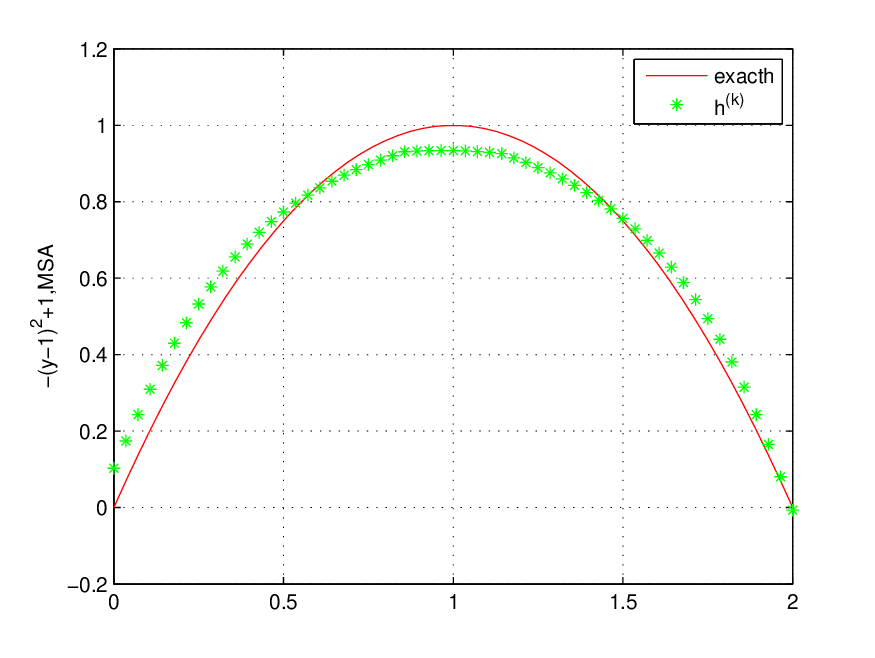} & \includegraphics[width=7cm]{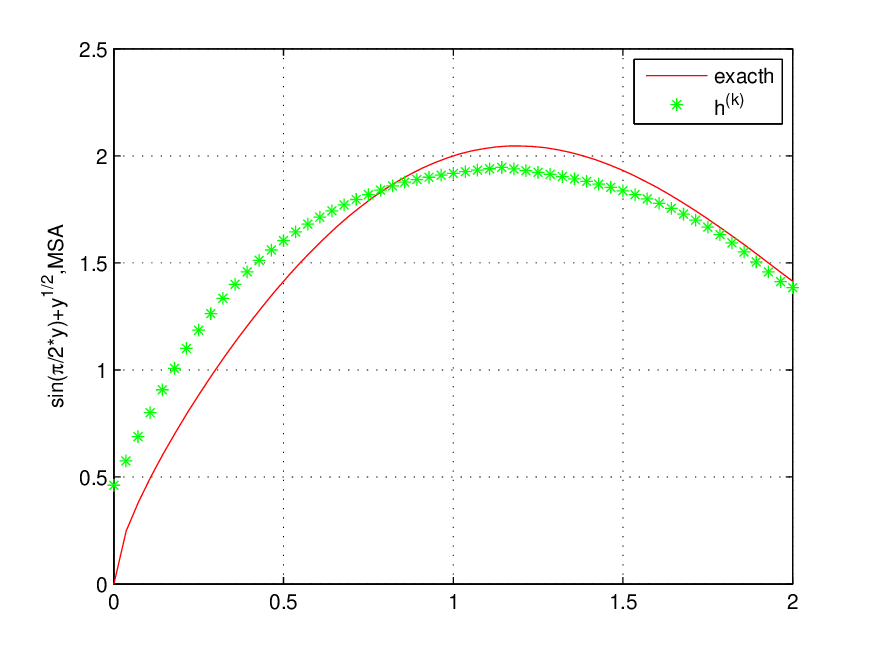} \\
%    (a)  Exact and reconstructed parameters & (b) exact data 
  \end{tabular}
  \caption{Exact and reconstructed fluxes for Examples 5.1 (left) and 
  5.2 (right)}\label{fig:ps}
\end{figure}

Figure \ref{fig:ps} (left) shows the exact parameter and the
numerically recovered parameter $h^{(k)}$, while Table \ref{s1} gives the number $k$ of iterations 
by Algorithms 3.1 (MSA) and 3.2 (ASA). 

%\begin{figure}[H]
%\caption{} \centerline{\includegraphics[width=300pt]{flux1.eps}}
%\label{f2}
%\end{figure}
%

\begin{table}
\caption{Iterative numbers $k$ of {\bf MSA} and {\bf ASA}}

~

\centering\begin{tabular}{cccccc}
 \hline Algorithm & N & M & $\beta$ & {error} & k
\\ \hline MSA & 14 & 28 & 0.0001 & 0.0597 & 8
\\ \hline $$ & 28 & 56 & 0.0001 & 0.0783 & 8
\\ \hline $$ & 56 & 112 & 0.0001 & 0.0907 & 8
\\ \hline\hline ASA  & 14 & 28 & 0.0001 & 0.0840 & 13
\\ \hline $$ & 28 & 56 & 0.0001 & 0.0978 & 13
\\ \hline $$ & 56 & 112 & 0.0001 & 0.0959 & 14
\\ \hline\end{tabular}
\label{s1}
\end{table}
% \hline Type & N & M & $\beta$ & error & k
%\\ \hline MSA   & 14 & 28 & 0.003 & 0.0428 & 31
%\\ \hline $$ & 28 & 56 & 0.003 & 0.0215 & 30
%\\ \hline $$ & 56 & 112 & 0.003 & 0.0264 & 30
%\\ \hline\hline MPA  & 14 & 28 & 0.01 & 0.0952 & 26
%\\ \hline $$ & 28 & 56 & 0.01 & 0.0931 & 27
%\\ \hline $$ & 56 & 112 & 0.01 & 0.0906 & 28

\begin{example}
We take the exact flux $h=\sin(\frac{\pi}{2}y)+y^{\frac{1}{2}}$ on $\Gamma_1$, the noise level $\delta=5\%$
and the constant initial guess $h^{(0)}=2$.
\end{example}

Figure \ref{fig:ps} (right) shows the exact parameter and the
numerically recovered parameter, while Table \ref{s2} gives the number of iterations 
by Algorithms 3.1 (MSA) and 3.2 (ASA).

\begin{table}
\caption{Iterative numbers $k$ of {\bf MSA} and {\bf ASA}}

~

\centering\begin{tabular}{cccccc}
 \hline Algorithm & N & M & $\beta$ & {error} & k
\\ \hline MSA & 14 & 28 & 0.0001 & 0.0827 & 9
\\ \hline $$ & 28 & 56 & 0.0001 & 0.0995 & 10
\\ \hline $$ & 56 & 112 & 0.0001 & 0.0996 & 11
\\ \hline\hline ASA & 14 & 28 & 0.0001 & 0.0981 & 14
\\ \hline $$ & 28 & 56 & 0.0001 & 0.0970 & 16
\\ \hline $$ & 56 & 112 & 0.0001 & 0.0921 & 15
\\ \hline\end{tabular}
\label{s2}
\end{table}
% \hline Type & N & M & $\beta$ & error & k
%\\ \hline MSA & 14 & 28 & 0.003 & 0.0455 & 32
%\\ \hline $$ & 28 & 56 & 0.003 & 0.0299 & 31
%\\ \hline $$ & 56 & 112 & 0.003 & 0.0272 & 30
%\\ \hline\hline MPA & 14 & 28 & 0.01 & 0.0965 & 29
%\\ \hline $$ & 28 & 56 & 0.01 & 0.0972 & 30
%\\ \hline $$ & 56 & 112 & 0.01 & 0.0936 & 31

%\begin{figure}[H]
%\caption{} \centerline{\includegraphics[width=300pt]{flux2.eps}} 
%\label{f3}
%\end{figure}
%

We can see from Figure \ref{fig:ps}  that the
numerical reconstructed fluxes, with a $5\%$ noise in the data, 
appear to be quite satisfactory, in view of 
the severe ill-posedness of the inverse flux problem. 
More importantly, we observe from Tables \ref{s1} and \ref{s2} that 
the convergence of the DD algorithms are nearly optimal with the refinement of 
the finite element mesh, i.e., 
the number of iterations grows very mildly with the mesh refinement.

\medskip
Next, we demonstrate three numerical examples of reconstructing the source strength
$f(\x)$ 
in the system (\ref{qi1f}),  with  $a(\x)=(x+y)/100$, $c(\x)=1$ in $\Om$
and $g(\x)=0$ on $\p\Om$. 
We start with a constant initial guess $f^{(0)}=0$ in $\Om$.

\begin{example}
We take  the exact source strength $f(x,y)=\sin(2\pi
x)\sin(2\pi y)$ and the noise level $\delta=1\%$.
\end{example}

\begin{figure}[H]
\centerline{\includegraphics[width=300pt]{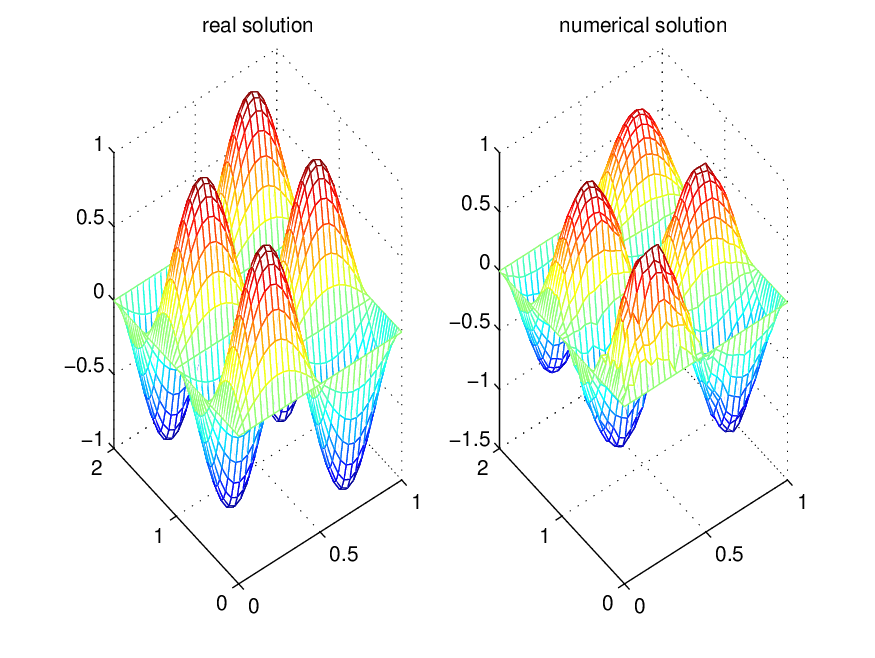}}
\caption{Exact and 
numerically recovered source strengths for Example 5.3} 
\label{aa1}
\end{figure}

Figure \ref{aa1} shows the exact and 
numerically recovered source strengths, while Table \ref{ss1} gives the number of iterations 
by Algorithms 2.1 (MSA) and 2.2 (ASA).

\begin{table}
\caption{Number of iterations of Algorithms 2.1 {(MSA)} and 2.2 {(ASA)}}

~

\centering\begin{tabular}{cccccc}
 \hline Algorithm & N & M & $\beta$ & error & k
\\ \hline MSA & 7 & 14 & 0.001 & 0.0900 & 10
\\ \hline $$  & 14 & 28 & 0.001 & 0.0971 & 14
\\ \hline $$ & 28 & 56 & 0.001 & 0.0998 & 14
\\ \hline $$ & 56 & 112 & 0.001 & 0.0979 & 15
\\ \hline\hline ASA & 7 & 14 & 0.001 & 0.0956 & 21
\\ \hline $$ & 14 & 28 & 0.001 & 0.0982 & 30
\\ \hline $$ & 28 & 56 & 0.001 & 0.0984 & 31
\\ \hline $$ & 56 & 112 & 0.001 & 0.0991 & 32
\\ \hline\end{tabular}
\label{ss1}
\end{table}
%\\ \hline\hline MPA & 7 & 14 & 0.001 & 0.0982 & 26
%\\ \hline $$ & 14 & 28 & 0.001 & 0.0970 & 28
%\\ \hline $$ & 28 & 56 & 0.001 & 0.0980 & 28
%\\ \hline $$ & 56 & 112 & 0.001 & 0.0983 & 29

\begin{example}
We take  the exact source strength $f=2\sin(2\pi x)y(y-1)(y-2)$
and the noise level $\delta=1\%$.
\end{example}

\begin{figure}[H]
\centerline{\includegraphics[width=300pt]{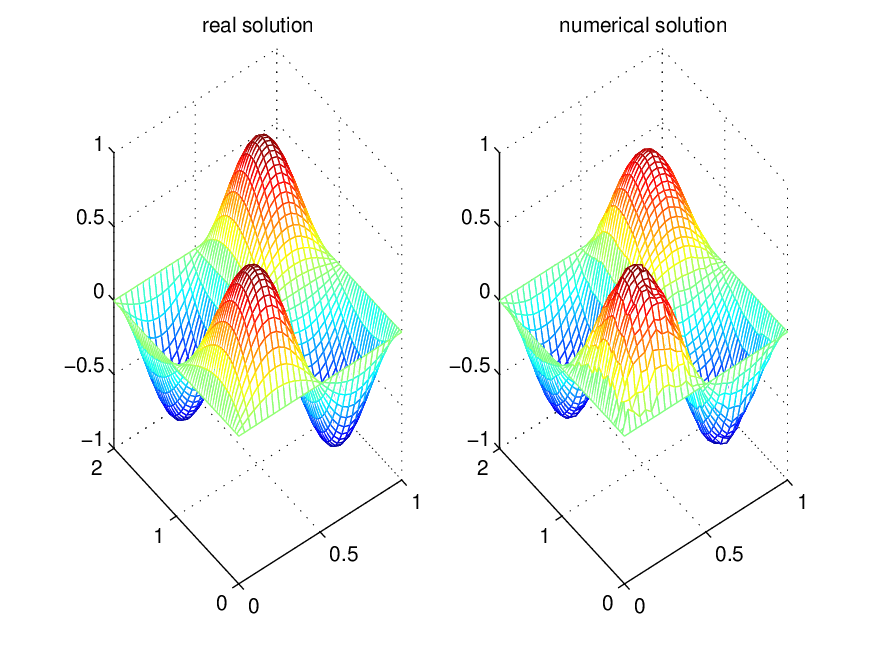}}
\caption{Exact and 
numerically recovered source strengths for Example 5.4} 
\label{aa2}
\end{figure}

Figure \ref{aa2} shows the exact and 
numerically recovered source strengths, while Table \ref{ss2} gives the number of iterations 
by Algorithms 2.1 (MSA) and 2.2 (ASA). 

%We take  $\lambda=\frac{1}{2}$ in MPA, the exact source strength
%$f=c_0*0.002\sin(2\pi x)y(y-1)(y-2)$, $c_0=1000$, the observed data
%$z^\delta=(1+\delta R)u(f)$, $\delta=0.01$ and the initial guess
%$\tilde f^{(0)}=0$ at all mesh points. Figure \ref{aa2} gives the exact source strength and the
%numerically recovered source strength $f^{(k)}$. One can see that the
%numerical reconstruction is very satisfactory. Table \ref{ss2} gives the
%iterative number $k$, when
%$\|f^{(k)}-f\|_{L^2(\Om)}/\|f\|_{L^2(\Om)}\leq 0.1$, we can see
%that $k$ is nearly stable with the refinement of meshes.

\begin{table}
\caption{Number of iterations for Algorithms 2.1 {(MSA)} and 2.2 {(ASA)}}

~

\centering\begin{tabular}{cccccc}
 \hline Algorithm & N & M & $\beta$ & error & k
\\ \hline MSA & 7 & 14 & 0.001 & 0.0933 & 6
\\ \hline $$  & 14 & 28 & 0.001 & 0.0991 & 7
\\ \hline $$ & 28 & 56 & 0.001 & 0.0895 & 8
\\ \hline $$ & 56 & 112 & 0.001 & 0.0970 & 8
\\ \hline\hline ASA & 7 & 14 & 0.001 & 0.0100 & 13
\\ \hline $$ & 14 & 28 & 0.001 & 0.0982 & 16
\\ \hline $$ & 28 & 56 & 0.001 & 0.0997 & 16
\\ \hline $$ & 56 & 112 & 0.001 & 0.0959 & 18
\\ \hline\end{tabular}
\label{ss2}
\end{table}
%\\ \hline\hline MPA & 7 & 14 & 0.001 & 0.0973 & 16
%\\ \hline $$ & 14 & 28 & 0.001 & 0.0991 & 14
%\\ \hline $$ & 28 & 56 & 0.001 & 0.0967 & 14
%\\ \hline $$ & 56 & 112 & 0.001 & 0.0985 & 14

\begin{example}
We take  the exact source strength $f(x,y)=10y\sin(2\pi y)x(x-1/2)(x-1)$
and the noise level $\delta=0.01$.
\end{example}

\begin{figure}[H]
\centerline{\includegraphics[width=300pt]{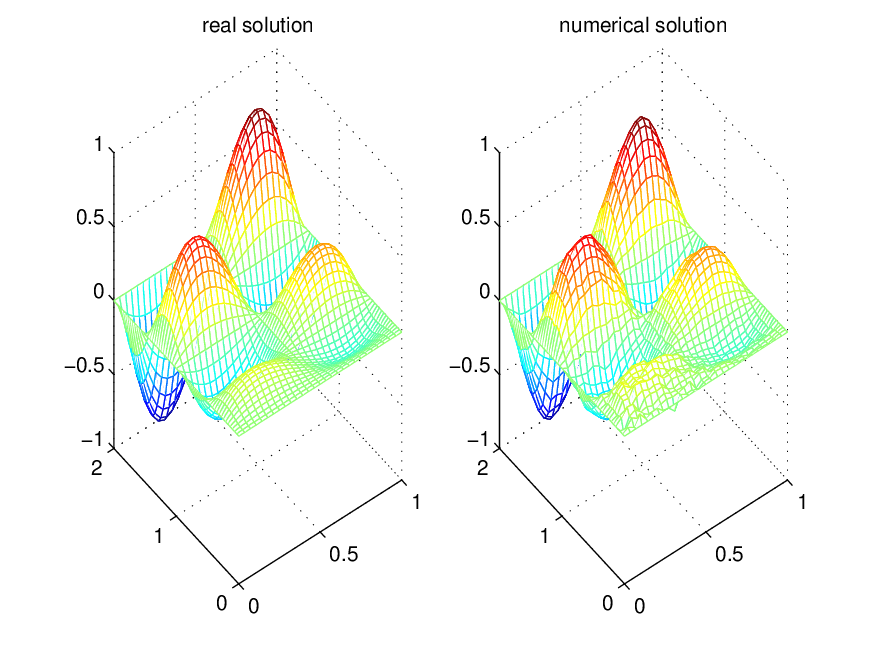}}
\caption{Exact and 
numerically recovered source strengths for Example 5.5} 
\label{aa3}
\end{figure}

Figure \ref{aa3} shows the exact and 
numerically recovered source strengths, while Table \ref{ss3} gives the number of iterations 
by Algorithms 2.1 (MSA) and 2.2 (ASA). 

\begin{table}
\caption{Number of iterations for Algorithms 2.1 {(MSA)} and 2.2 {(ASA)}}

~

\centering\begin{tabular}{cccccc}
 \hline Type & N & M & $\beta$ & error & k
\\ \hline MSA & 7 & 14 & 0.001 & 0.0989 & 14
\\ \hline $$  & 14 & 28 & 0.001 & 0.0976 & 23
\\ \hline $$ & 28 & 56 & 0.001 & 0.0975 & 24
\\ \hline $$ & 56 & 112 & 0.001 & 0.0980& 26
\\ \hline\hline ASA & 7 & 14 & 0.001 & 0.0977 & 30
\\ \hline $$ & 14 & 28 & 0.001 & 0.0989 & 47
\\ \hline $$ & 28 & 56 & 0.001 & 0.0997 & 48
\\ \hline $$ & 56 & 112 & 0.001 & 0.0999 & 52
\\ \hline\end{tabular}
\label{ss3}
\end{table}

%\begin{remark}
%We can see from Figures \ref{aa1}-\ref{aa3} that the
%numerical identification for source strength are very satisfactory.
%Tables \ref{ss1}-\ref{ss3} show that the
%iterative number $k$ when algorithms stop grow very slowly with the refinement of meshes.
%\end{remark}
%
We can see from Figures \ref{aa1}-\ref{aa3} that the
numerical reconstructed source strengths, with a $1\%$ noise in the data, 
appear to be quite satisfactory, in view of 
the severe ill-posedness of the inverse source problem and the complicated profiles 
of the exact source strengths, especially in Example 5.3 where the source strength 
oscillates frequently between 8 peaks and valleys. 
More importantly, we observe from Tables \ref{ss1}-\ref{ss3} that 
the convergence of the DD algorithms are nearly optimal with the refinement of 
the finite element mesh, i.e., 
the number of iterations grows only mildly with the mesh refinement.

\medskip
Finally, we present three numerical examples for the reconstructions of 
the initial temperature in the heat conductive system (\ref{ya1}),  by two DD 
algorithms, namely Algorithms 4.1 and 4.2 proposed in Section \ref{sec:initial}. 
In our experiments, we take 
$a(\x)=1$, $f(\x,t)=0$, the terminal time $T=4$, with the constant initial
guess $\varphi^{(0)}=0$.

\begin{example}
We take  the exact initial temperature
$\varphi=\sin(2\pi x)\sin(2\pi y)$ and the noise level $\delta=2\%$.
\end{example}

\begin{figure}[H]
\centerline{\includegraphics[width=300pt]{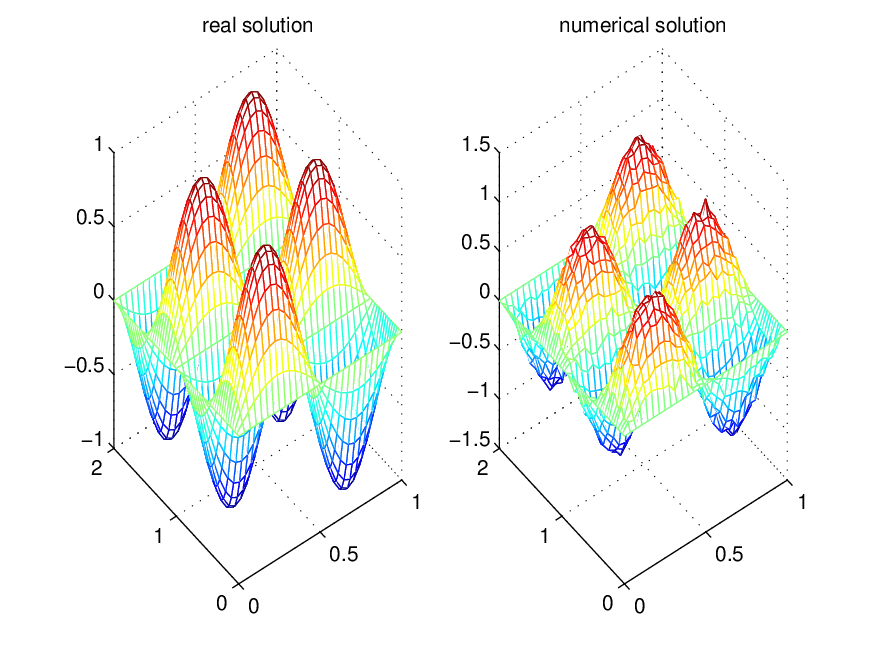}}
\caption{Exact and 
numerically recovered initial temperatures for Example 5.6}
\label{bb1}
\end{figure}

Figure \ref{bb1} shows the exact and 
numerically recovered initial temperatures, while Table \ref{t1} gives the number of iterations 
by Algorithms 4.1 (MSA) and 4.2 (ASA). 

\begin{table}
\caption{Number of iterations for Algorithms 4.1 {(MSA)} and 4.2 {(ASA)}}

~

\centering\begin{tabular}{cccccc}
 \hline Algorithm & N & M & $\beta$ & {error} & k
\\ \hline MSA & 7 & 14 & 0.00005 & 0.0521 & 10
\\ \hline $$  & 14 & 28 & 0.00005 & 0.0760 & 11
\\ \hline $$ & 28 & 56 & 0.00005 & 0.0998 & 12
%\\ \hline $$ &56 ¡¡¡¡¡¡¡¡& 112 & 0.00005 & - & -
\\ \hline\hline ASA & 7 & 14 & 0.00005 & 0.0856 & 16
\\ \hline $$ & 14 & 28 & 0.00005 & 0.0950 & 20
\\ \hline $$ & 28 & 56 & 0.00005 & 0.0973 & 25
%\\ \hline $$ &56 ¡¡¡¡¡¡¡¡& 112 & 0.00005 & - & -
\\ \hline\end{tabular}
\label{t1}
\end{table}

\begin{example}
We take  the exact initial temperature
$\varphi=2\sin(2\pi x)y(y-1)(y-2)$ and the noise level $\delta=1\%$.
\end{example}

\begin{figure}[H]
\centerline{\includegraphics[width=300pt]{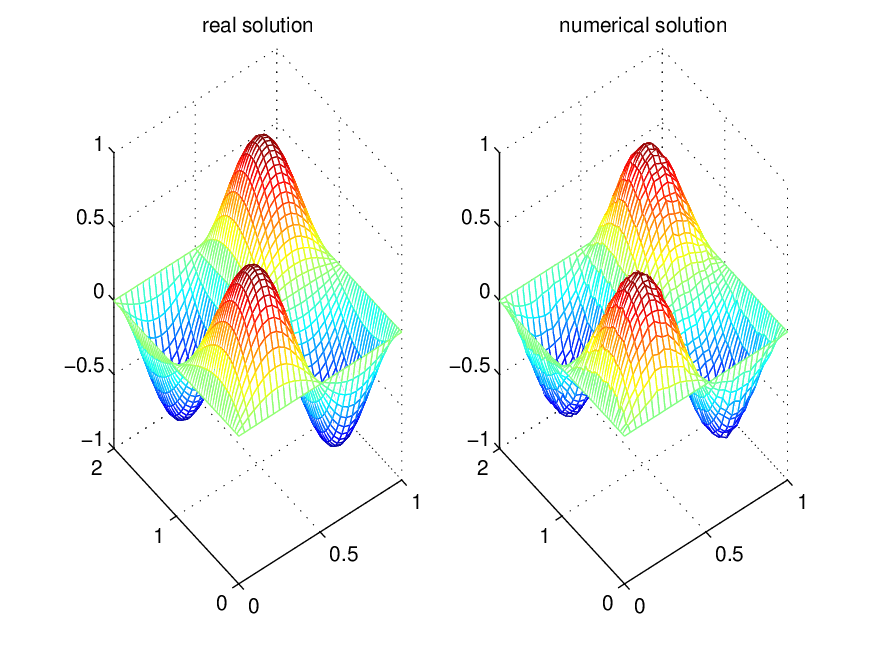}}
\caption{Exact and 
numerically recovered initial temperatures for Example 5.7}
\label{bb2}
\end{figure}

Figure \ref{bb2} shows the exact and 
numerically recovered initial temperatures, while Table \ref{t2} gives the number of iterations 
by Algorithms 4.1 (MSA) and 4.2 (ASA). 

\begin{table}
\caption{Number of iterations of Algorithms 4.1 (MSA) and 4.2 (ASA)}

~

\centering\begin{tabular}{cccccc}
 \hline Algorithm & N & M & $\beta$ & {error} & k
\\ \hline MSA & 7 & 14 & 0.00005 & 0.0943 & 22
\\ \hline $$  & 14 & 28 & 0.00005 & 0.0931 & 26
\\ \hline $$ & 28 & 56 & 0.00005 & 0.0995 & 29
%\\ \hline $$ &56 ¡¡¡¡¡¡¡¡& 112 & 0.00005 & - & -
\\ \hline\hline ASA & 7 & 14 & 0.00005 & 0.0974 & 45
\\ \hline $$ & 14 & 28 & 0.00005 & 0.0997 & 52
\\ \hline $$ & 28 & 56 & 0.00005 & 0.0971 & 58
%\\ \hline $$ &56 ¡¡¡¡¡¡¡¡& 112 & 0.00005 & - & -
\\ \hline\end{tabular}
\label{t2}
\end{table}

\begin{example}
We take  the exact initial temperature
$\varphi=10y\sin(2\pi y)x(x-1/2)(x-1)$ and the noise level $\delta=2\%$.
\end{example}

\begin{figure}[H]
\centerline{\includegraphics[width=300pt]{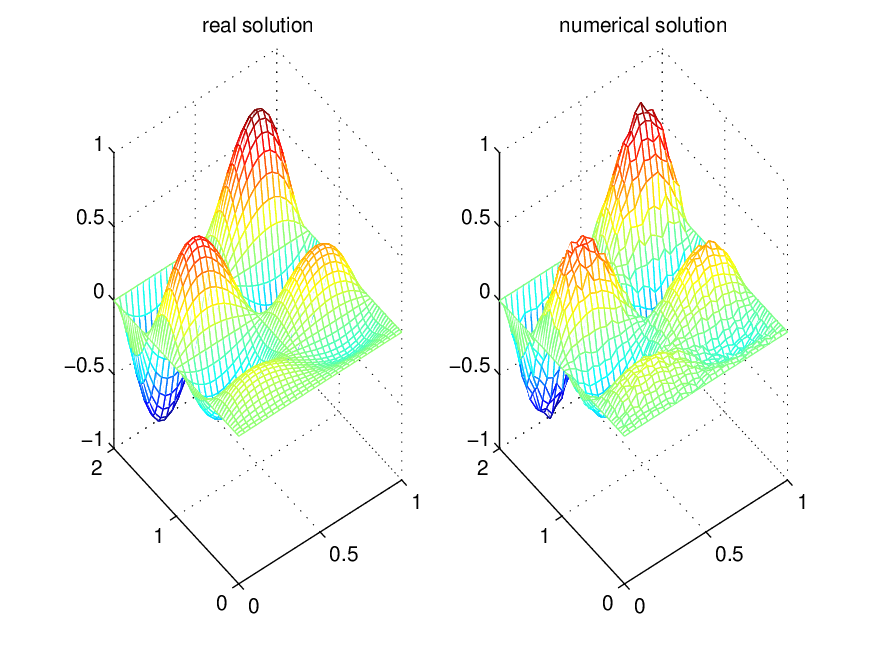}}
\caption{Exact and 
numerically recovered initial temperatures for Example 5.8}
\label{bb3}
\end{figure}

Figure \ref{bb3} shows the exact and 
numerically recovered initial temperatures, while Table \ref{t3} gives the number of iterations 
by Algorithms 4.1 (MSA) and 4.2 (ASA). 

\begin{table}
\caption{Number of iterations for Algorithms 4.1 (MSA) and 4.2 (ASA)}

~

\centering\begin{tabular}{cccccc}
 \hline Algorithm & N & M & $\beta$ & {error} & k
\\ \hline MSA & 7 & 14 & 0.00005 & 0.0933 & 8
\\ \hline $$  & 14 & 28 & 0.00005 & 0.0956 & 10
\\ \hline $$ & 28 & 56 & 0.00005 & 0.0989 & 15
%\\ \hline $$ &56 ¡¡¡¡¡¡¡¡& 112 & 0.00005 & -& -
\\ \hline\hline ASA & 7 & 14 & 0.00005 & 0.0988 & 17
\\ \hline $$ & 14 & 28 & 0.00005 & 0.0970 & 22
\\ \hline $$ & 28 & 56 & 0.00005 & 0.0993 & 30
%\\ \hline $$ &56 ¡¡¡¡¡¡¡¡& 112 & 0.00005 & - & -
\\ \hline\end{tabular}
\label{t3}
\end{table}

We can see from Figures \ref{bb1}-\ref{bb3} that the
numerical reconstructed initial temperatures, with a $1\%$ noise in the data, 
appear to be quite satisfactory, in view of 
the severe ill-posedness of the inverse initial temperature problem and the complicated profiles 
of the exact initial temperatures, especially in Example 5.6 where the initial temperature
oscillates frequently between 8 peaks and valleys. 
More importantly, we observe from Tables \ref{t1}-\ref{t3}  that 
the convergence of the DD algorithms are nearly optimal with the refinement of 
the finite element mesh, i.e., 
the number of iterations grows still mildly with the mesh refinement.
But compared with the numerical results for the source strengths and fluxes, 
we can see that the performance of the reconstructions for the initial temperatures 
are less satisfactory in terms of the mesh refinement.

\section{Concluding remarks}\l{sec:conclusion}
\setcounter{equation}{0}
We have proposed several overlapping domain decomposition algorithms 
for solving some representative linear inverse problems,
including the identification of the fluxes, the source intensity and the initial temperature
in second order elliptic and parabolic systems.
The algorithms are constructed in a way that only
small sub-minimizations are needed to solve on the subdomains of the original global domain at each iteration.
And it is important to observe from many numerical examples 
that the convergence of the DD algorithms are nearly optimal with the refinement of 
the finite element mesh, i.e.,  the number of iterations grows only mildly with the mesh refinement.

Our future work includes the extension of the proposed overlapping domain decomposition algorithms 
to nonlinear inverse problems, such as the constructions of the diffusivity coefficient,
the radiative coefficient and Robin coefficient in elliptic and parabolic systems.

\end{document}